\documentclass[11pt,oneside]{amsart}

\usepackage[T1]{fontenc}
\usepackage[utf8]{inputenc}

\usepackage{amsmath}
\usepackage{amssymb}

\usepackage{enumitem}
\usepackage{graphicx}

\usepackage[colorlinks]{hyperref}

\newtheorem{thm}{Theorem}[section]

\newtheorem{cor}[thm]{Corollary}

\newtheorem{exmpl}[thm]{Example}
\newtheorem{lem}[thm]{Lemma}
\newtheorem{prop}[thm]{Proposition}

\newtheorem{rem}[thm]{Remark}

\newtheorem{ques}[thm]{Question}

\newtheorem{fact}[thm]{Fact}

\def\L{\mathcal{L}}

\DeclareMathOperator{\Span}{Span}
\DeclareMathOperator{\Image}{Im}
\DeclareMathOperator{\Rauz}{Rauz}

\newcommand*\W{\mathcal{W}}

\title[Polynomial complexity, random graphs, and strict ergodicity]{Polynomial complexity, random graphs, and strict ergodicity}

\author{Be'eri Greenfeld} \address{Department of Mathematics and Statistics, Hunter College and CUNY Graduate Center, 695 Park Avenue, New York NY 10065, USA} \email{beeri.greenfeld@hunter.cuny.edu}

\author{Paul-Henry Leemann}
\email{paul-henry.leemann@protonmail.com}

\author{Carlos Gustavo Moreira} \address{SUSTech International Center for Mathematics, Shenzhen, Guangdong,
People’s Republic of China;
IMPA, Estrada Dona Castorina 110, 22460-320, Rio de Janeiro, Brazil} \email{gugu@impa.br}

\author{Efim Zelmanov} \address{SUSTech International Center for Mathematics, Shenzhen, Guangdong,
People’s Republic of China} \email{efim.zelmanov@gmail.com}

\begin{document}

\begin{abstract}
We construct infinite words $w$ whose complexity functions are  quadratically bounded, but whose discrete derivatives are not linearly bounded, resolving a question of Cassaigne.

Moreover, we construct infinite words of polynomially bounded complexity, and subshifts of quasi-linear complexity, whose discrete derivatives do not decay compared to their complexity functions, thus exhibiting a typical exponentially-growing behavior.
As an application, we give the first examples of subshifts of polynomial complexity whose Rauzy graphs do not have Benjamini-Schramm limits.

Finally, we realize (up to an inevitable linear error factor) every function satisfying the necessary conditions, as the complexity function of a minimal, uniquely ergodic, cube-free subshift, strengthening previous realization results in this direction, and derive consequences for the growth of simple algebras and just-infinite semigroups satisfying nontrivial laws.
\end{abstract}

\maketitle

\section{Introduction}

Let $\Sigma$ be a finite alphabet and let $$w=w_1w_2\cdots \in \Sigma^{\mathbb{N}}$$
be an infinite word over $\Sigma$. The \emph{complexity function} of $w$, denoted $p_w(n)$, counts the number of subwords (or `factors') of $w$ of length $n$. This is a natural and useful measure of `how complicated' a given infinite word is. More generally, a \emph{subshift} is a non-empty closed, shift-invariant subspace of $\Sigma^{\mathbb{N}}$. The closure of the shift orbit of any infinite word is a subshift; complexity functions of subshifts are defined similarly. 
A typical theme in symbolic dynamics is that complexity functions reflect, and to a certain extent control, the dynamics of a given infinite word or subshift.

A fundamental problem in combinatorial symbolic dynamics is to understand how flexible the complexity functions of infinite words and subshifts can be, and how flexible subshifts of restricted complexity can be (a partial list of references includes \cite{All94, BBL, BST19, Cas2, Cas, Cas3,CP23, CPsubmitted,CK15,CK19,DLR13,Fer99,Fer96,GMZ,MM1,MM2,  MH38,MH} and references therein). In \cite{GMZ}, the first, third and fourth named authors characterized these functions up to bi-Lipschitz equivalence.

Considering the fine asymptotics of complexity functions beyond bi-Lipschitz equivalence, additional interesting regularity patterns appear. Given a function $f\colon \mathbb{N} \to \mathbb{N}$ we denote its discrete derivative by $f'(n)=f(n)-f(n-1)$. 
Ferenczi conjectured that if $p_w=O(n)$ then $p_w'(n)=O(1)$. This was confirmed by Cassaigne \cite{Cas3}, who raised a more general question:
\begin{ques} \label{ques:Cas}
If $p_w(n)=O(n^\alpha)$, must $p_w'(n)=O(n^{\alpha-1})$?
\end{ques}
Thus, Cassaigne's theorem from \cite{Cas3} solves this question in the affirmative for $\alpha=1$. Moreover, in the same paper, Cassaigne proved that if $p_w(n)=O(n^\alpha)$ then $p_w'(n)=O(n^{3\alpha-3})$.

Our first main result resolves Question \ref{ques:Cas} in the negative.

\begin{thm} \label{thm:quadratic}
For every $a>1$, there exists an infinite word $w$ such that $$p_w(n)=O(n^a)\ \ \ \text{but}\ \ \ p_w'(n)\neq O(n^{a-1}).$$
\end{thm}
A question related to Question \ref{ques:Cas} is 
\begin{ques} \label{question}
Let $w$ be an infinite word of polynomially bounded (or: subexponential) complexity. Is $\frac{p_w'(n)}{p_w(n)}\xrightarrow{n\rightarrow \infty} 0$?
\end{ques}
A similar question can be asked for arbitrary subshifts. Notice that the asymptotic behavior implied in the above question is indeed typical for subexponential functions. 
The motivation behind Question \ref{question} comes from the study of Benjamini-Schramm limits of \emph{Rauzy graphs}. 
The Rauzy graph $\text{Rauz}_\mathcal{X}(n)$ ($n=1,2,\dots$) associated with a subshift $\mathcal{X}$ is the finite directed graph whose vertices are the length-$n$ factors of $\mathcal{X}$, and whose edges correspond to length-$n+1$ factors: the vertices $xu,uy$ ($x,y$ are letters and $|u|=n-1$) are connected by an edge if $xuy$ is again a factor of $\mathcal{X}$. Thus, Rauzy graphs give finite approximations of subshifts.
By choosing a root uniformly at random, each finite graph induces a finitely supported measure on the space of rooted graphs. If the sequence of measures associated with a sequence of finite graphs has a weak limit, thought of as a random graph, it is called the \emph{Benjamini-Schramm} limit of the given sequence of graphs.
Benjamini-Schramm limits of Rauzy graphs associated with subshifts have been extensively studied, see e.g. \cite{GLN,L2016,LNSV, LNV}.
It was proved in \cite{LNSV} that if $p_\mathcal{X}(n)$ is subexponential and the language of \(\mathcal{X}\) is almost-prolongable
then $$\lim_{n\to \infty}\Rauz_{\mathcal{X}}(n)\ \text{exists}\ \iff \ \lim_{n \to \infty}\frac{p_\mathcal{X}'(n)}{p_\mathcal{X}(n)}=0.$$
A language is \emph{almost prolongable} if the proportion of words of length \(n\) that are both left and right prolongable goes to \(1\) when \(n\) goes to infinity. Languages of two-sided subshifts or of one-sided subshifts of unbounded complexity with a finite subset whose shift orbit is dense are almost prolongable \cite{LNSV}.
Question \ref{question} can be now interpreted as:
\begin{ques} \label{question2}
Let $\mathcal{X}$ be a subshift of polynomially bounded (or: subexponential) complexity and whose language is almost prolongable.
Does $\lim_{n\to \infty} \Rauz_\mathcal{X}(n)$ exist?
\end{ques}
Affirmative answers to Questions \ref{question}, \ref{question2} have been proved in some important cases: if the complexity is linear \cite{Cas3}; more generally, if the complexity is $O(n^{\alpha})$ for any $\alpha<\frac{4}{3}$;
and for regular formal languages (that is, those accepted by finite automata) of subexponential complexity.

Our second main result gives the first counterexamples to Question \ref{question}.
\begin{thm} \label{thm:p'>cp}
    There exists a recurrent infinite word $w$ whose complexity function satisfies $$p_w(n) \leq n^{3+o(1)}$$
    but for some universal constant $\delta>0$,  $$p_w'(n) > \delta \cdot p_w(n)\ \text{infinitely often}.$$
\end{thm}
Consequently, we derive the first counterexamples to Question \ref{question2}:

\begin{cor}
There exists a transitive (irreducible) subshift $\mathcal{X}$ of polynomial complexity for which $\lim_{n\to \infty} \Rauz_\mathcal{X}(n)$ does not exist.
\end{cor}
\noindent (Notice that Theorem \ref{thm:p'>cp} does not imply a negative answer to Question \ref{ques:Cas}, which indeed requires a different approach.)
If the set of indices $n$ for which $p_w'(n)>c\cdot p_w(n)$ for some constant $c>0$ is large enough (has a positive lower density) then $p_w(n)$ grows exponentially (see Remark \ref{remark}). In \cite{GMZ2}, the first, third and fourth named authors gave examples of infinite words $w$ with polynomially bounded complexity for which $p_w'(n)\neq O(p_w(n)/n^\varepsilon)$. Theorem \ref{thm:p'>cp} significantly improves this.

A phenomenon as in Theorem \ref{thm:p'>cp} can occur even beyond bi-Lipschitz equivalence, namely, when the argument of $p_w$ can be scaled:

\begin{prop} \label{prop:p'>pc}
There exists a subshift $\mathcal{X}\subseteq \{0,1\}^{\mathbb{N}}$ such that $p_\mathcal{X}(n)$ is quasi-linear (in fact, $O(n \log^2 n)$) but for every $C\geq 1$ there exists $c>0$ such that $$p_\mathcal{X}'(n) \geq c \cdot p_\mathcal{X}(Cn)$$ infinitely often.
\end{prop}
(The subshifts of quasi-linear complexity constructed above are one-sided and do not have a finite subset whose shift orbits are dense. This is inevitable: indeed, otherwise, the limit exists by \cite[Corollary~1.2]{LNSV}.)

\medskip

We then turn to construct subshifts and infinite words with prescribed complexity and with extreme combinatorial and dynamical properties. 
In \cite{GMZ,GMZ2}, the first, third and fourth named authors constructed, for any non-decreasing submultiplicative unbounded function $f\colon \mathbb{N}\rightarrow \mathbb{N}$, a minimal (even strictly ergodic) subshift $\mathcal{X}$ whose complexity function `approximates' $f$, that is, $f(n) \preceq p_\mathcal{X}(n) \preceq nf(n)$. Borrowed from large-scale geometry and from geometric group theory, we write $f(n)\preceq g(n)$ if, for some constants $C,D>0$, for all $n\gg 1$ we have $f(n)\leq Cg(Dn)$; and $f(n)\sim g(n)$ are \emph{bi-Lipschitz equivalent} if both $f(n)\preceq g(n)$ and $g(n)\preceq f(n)$.

We now close this `approximation gap' by proving that, in fact, we can take $p_\mathcal{X}(n) \sim nf(n)$ (the linear error factor is necessary in this formulation, as complexity functions are either bounded or at least linear). This is achieved by a certain `padding' process, building on our previous `subshift factory' and inspired by an earlier idea of Currie--Rampersad \cite{CR10} of interleaving infinite words with the Thue-Morse sequence.

\begin{thm} \label{thm:ergodic}
    Let $f\colon \mathbb{N}\rightarrow \mathbb{N}$ be a non-decreasing, unbounded, submultiplicative function. Then there exists a minimal, uniquely ergodic, cube-free subshift $\mathcal{X}$ such that $$p_\mathcal{X}(n) \sim nf(n).$$
\end{thm}

We conclude with applications to the growth of algebraic structures. 

\begin{cor} \label{cor:simple}
    Let $f\colon \mathbb{N}\rightarrow \mathbb{N}$ be a non-decreasing, submultiplicative function. Then there exists a finitely generated simple associative algebra $A$ such that $$\gamma_A(n)\sim n^2 f(n).$$
\end{cor}

\begin{cor}
    Let $f\colon \mathbb{N}\rightarrow \mathbb{N}$ be a non-decreasing, submultiplicative function. Then there exists a finitely generated just-infinite semigroup with zero $S$ satisfying a non-trivial semigroup law, namely, $X^3=0$, of growth $$\gamma_S(n)\sim n^2 f(n).$$
\end{cor}

\paragraph{Acknowledgement}
The second-named author was supported by the RDF-23-01-045 grant from XJTLU. He would also like to extend his warmest thanks to SUSTech for their hospitality. 
The fourth-named author was partially supported by the NSFC Grant 12350710787 and the Guangdong Program 2023JC10X085.

\section{Preliminaries on combinatorics on words}

Let $\Sigma$ be a finite alphabet. Let $\mathcal{L}\subseteq \bigcup_{n=0}^{\infty} \Sigma^n$ be a formal language. We denote $\mathcal{L}_n:=\mathcal{L}\cap \Sigma^n$, so $p_\mathcal{L}(n) = |\mathcal{L}_n|$. 
A formal language is hereditary if it is closed under taking subwords. A hereditary language $\mathcal{L}$ is right (respectively, left) prolongable if for every $u\in \mathcal{L}$ there exists some $a\in \Sigma$ such that $ua\in \mathcal{L}$ (resp. $au\in \mathcal{L}$), and prolongable if it is right and left prolongable. 

Given $m \leq n$ we let $\text{pre}^n_m\colon \Sigma^n\rightarrow \Sigma^{m}$ (respectively, $\text{suf}^n_m\colon \Sigma^n\rightarrow \Sigma^{m}$) denote the natural prefix (resp., suffix) map. If $n,m$ are obvious from context, we omit the superscript and subscript. Thus a hereditary language is right (resp., left) prolongable if and only if $\text{pre}\colon \mathcal{L}_n\rightarrow \mathcal{L}_{n-1}$ (resp., $\text{suf}\colon \mathcal{L}_n\rightarrow \mathcal{L}_{n-1}$) is surjective for all $n$. Equivalently, $\mathcal{L}$ coincides with the hereditary language of a right (resp., left) one-sided subshift; and $\L$ is prolongable if and only if it coincides with the hereditary language of a bi-infinite subshift $\mathcal{X}\subseteq \Sigma^\mathbb{Z}$. A subshift is transitive if and only if it is the closure of the shift orbit of a single infinite word. The complexity function of any subshift is non-decreasing.

An infinite word $w$ (one-sided or two-sided) is \emph{recurrent} if every subword of it appears infinitely often (on both sides, in the two-sided case). A subshift $\mathcal{X}$ is \emph{irreducible} if, whenever $u,v$ are factors of $\mathcal{X}$, there exists some $w$ such that $uwv$ is again a factor of $\mathcal{X}$.
A subshift is irreducible if and only if it is generated by a recurrent infinite word. An infinite word $w$ is uniformly recurrent if for every subword $u$ of $w$ there exists some $N_u$ such that every subword of $w$ of length $N_u$ contains an occurrence of $u$. A word is \emph{uniformly recurrent} if and only if it generates a minimal subshift, namely, a subshift with no proper subshifts. The classes of hereditary languages associated with one-sided recurrent words and with two-sided recurrent words, coincide.

Given two subsets $S,T\subseteq \bigcup_{n=0}^{\infty} \Sigma^n$, we let $S\cdot T = \{st\ |\ s\in S,\ t\in T\}$.

As mentioned in the introduction, given two non-decreasing functions $f,g\colon \mathbb{N}\rightarrow \mathbb{N}$, we write $f(n)\preceq g(n)$ if, for some constants $C,D>0$, for all $n\gg 1$ we have $f(n)\leq Cg(Dn)$; and $f(n)\sim g(n)$ are \emph{bi-Lipschitz equivalent} if both $f(n)\preceq g(n)$ and $g(n)\preceq f(n)$. In other words, if there exists some constant $c>0$ such that $\frac{1}{c}f(\frac{1}{c}n)\leq g(n) \leq cf(cn)$. Up to bi-Lipschitz equivalence, all exponential functions coincide. Bi-Lipschitz equivalence can be determined on sets with bounded ratios between successive elements.

\begin{lem}
Let $f,g\colon \mathbb{N}\rightarrow \mathbb{N}$ be non-decreasing functions. Let $1\leq s_1<s_2<\dots$ be an infinite sequence in which, for some constant $K>0$, $1<s_{i+1}/s_i<K$ for all $i$. Suppose that for some constants $C,D>0$, we have $f(s_i)\leq Cg(Ds_i)$ for all $i$. Then $f(n)\preceq g(n)$.
\end{lem}
\begin{proof}
Given $n\geq s_1$, there exists some $i$ such that $s_i \leq n <s_{i+1}$, so
$$f(n)\leq f(s_{i+1}) \leq C g(D s_{i+1}) \leq C g(D K n)$$
since $\frac{s_{i+1}}{n}\leq \frac{s_{i+1}}{s_i} \leq K$.
\end{proof}

We say that a non-decreasing function $f\colon \mathbb{N}\rightarrow \mathbb{N}$ is subexponential if for every $\lambda>1$, we have $\frac{f(n)}{\lambda^n}\xrightarrow{n\rightarrow \infty} 0$; equivalently, $f(n)^{1/n}\xrightarrow{n\rightarrow \infty} 1$. 

The following is a well-known consequence of Fekete's Lemma; we bring it here for completeness.

\begin{lem} \label{lem:exponential always}
Let $f\colon \mathbb{N}\rightarrow \mathbb{N}$ be a non-decreasing, submultiplicative function. If there exists a subexponential function $g\colon \mathbb{N}\rightarrow \mathbb{N}$ such that $f(n) \leq g(n)$ infinitely often then $f(n)$ is subexponential.
\end{lem}
\begin{proof}
    By Fekete's Lemma, the limit $\lambda := \lim_{n\rightarrow \infty} f(n)^{1/n}$ exists. Since $f(n_k)^{1/n_k}\leq g(n_k)^{1/n_k}\xrightarrow{k\rightarrow \infty} 1$, it follows that $\lambda=1$, so $f$ is subexponential.
\end{proof}

\begin{rem} \label{remark}
    Let $f\colon \mathbb{N}\rightarrow \mathbb{N}$ be a non-decreasing submultiplicative function (e.g. the complexity function of any subshift). Suppose that there exists some constant $c>0$ such that $I:=\{n\in \mathbb{N}|f'(n)>cf(n)\}$ has a positive lower density $\delta>0$, then $f(n)$ grows exponentially. Indeed, $I\subseteq \{n\in \mathbb{N}|\frac{f(n)}{f(n-1)}>1+c\}$ and so $$f(N)=f(1)\prod_{k=2}^{N} \frac{f(k)}{f(k-1)}\geq (1+c)^{|I\cap [2,N]|} \geq (1+c)^{0.99\delta N},$$ (the last inequality holding for $N\gg 1$) which grows exponentially.
\end{rem}

\begin{lem} \label{lem:derivative}
Let $\mathcal{L}$ be a right prolongable hereditary language over a finite alphabet (namely, the hereditary language of an $\mathbb{N}$-subshift). Then
$$
p_\mathcal{L}(n)-p_\mathcal{L}(n-d) \geq \Big|\Big\{\substack{\text{Words of length}\ n-d\ \text{with at least two} \\ \text{right prolongations by}\ d\ \text{letters}}\Big\}\Big|.
$$
In particular,
$$
p_\mathcal{L}'(n) \geq \Big|\Big\{\substack{\text{Words of length}\ n-1\ \text{with at least two} \\ \text{right prolongations by a single letter}}\Big\}\Big|.
$$
\end{lem}
\begin{proof}
Given $m,d$, let $\mathcal{L}_m^{(+d)} \subseteq \mathcal{L}_m$ be the set of length-$m$ words in $\mathcal{L}$ that have at least two distinct length-$d$ right prolongations.
This means that, for the prefix map $\text{pre}^{n}_{n-d}\colon \mathcal{L}_n \rightarrow \mathcal{L}_{n-d}$ we have
\begin{align*}
    & |(\text{pre}^n_{n-d})^{-1}(\mathcal{L}_{n-d}^{(+d)})| \geq 2\cdot |\mathcal{L}_{n-d}^{(+d)}| \ \text{and} \\ & |(\text{pre}^n_{n-d})^{-1}(\mathcal{L}_{n-d} \setminus \mathcal{L}_{n-d}^{(+d)})| =|\mathcal{L}_{n-d} \setminus \mathcal{L}_{n-d}^{(+d)}| = p_\mathcal{L}(n-d) - |\mathcal{L}_{n-d}^{(+d)}|
\end{align*}
so \begin{align*}
    p_\mathcal{L}(n) = |\mathcal{L}_n| & = |(\text{pre}^{n}_{n-d})^{-1}(\mathcal{L}_{n-d})| \\ & = |(\text{pre}^n_{n-d})^{-1}(\mathcal{L}_{n-d}^{(+d)})| + |(\text{pre}^n_{n-d})^{-1}(\mathcal{L}_{n-d} \setminus \mathcal{L}_{n-d}^{(+d)})| \\ & \geq 2\cdot |\mathcal{L}_{n-d}^{(+d)}| + \left( p_\mathcal{L}(n-d) - |\mathcal{L}_{n-d}^{(+d)}| \right)=p_\mathcal{L}(n-d) + |\mathcal{L}_{n-d}^{(+d)}|
\end{align*}
hence $$ p_\mathcal{L}(n) - p_\mathcal{L}(n-d) \geq |\mathcal{L}_{n-d}^{(+d)}|, $$
as claimed.
\end{proof}

\noindent \emph{Conventions.} We occasionally abbreviate the word $\underbrace{a\cdots a}_{k}$ by $a^k$ for brevity when it improves the readability of the statement. For functions $f,g\colon \mathbb{N}\rightarrow \mathbb{N}$, we write $f(n)=O(g(n))$ if $f(n)\leq Cg(n)$ for some constant $C>0$, and $f(n)=o(g(n))$ if $f(n)/g(n)\xrightarrow{n\rightarrow \infty} 0$.

\section{Words of quadratic complexity and super-linear derivative}

Recall that a (binary) \emph{De Bruijn sequence} of order $n$ is a binary sequence in which, when written cyclically, every length-$n$ binary sequence appears exactly once. The length of such a sequence is $2^n$. 
For example, $$00010111$$ is a binary De Bruijn sequence of order $3$, and
$$0000100110101111$$
is a binary De Bruijn sequence of order $4$. 
Let $$t_n \colon \mathbb{Z}/2^n\mathbb{Z} \longrightarrow \{0,1\}$$ be a binary De Bruijn sequence of order $n$; we may assume that $t_n(0)=0$.

\medskip

We work over the alphabet $\Sigma = \{0,1,2\}$; we let $\tau$ be the substitution given by $\tau(0)=0,\ \tau(1)=2,\ \tau(2)=1$.
We fix the following two sequences of positive integers, $\{D_n\}_{n=1}^{\infty}$ (a sequence of `sizes') and $\{k_n\}_{n=1}^{\infty}$ (a sequence of `repetition orders'). We take $D_1=k_1=1$ and iteratively set $D_{n+1}=2^{k_n+1}D_n$; the sequence of repetition orders $\{k_n\}_{n=1}^{\infty}$ will be determined in the sequel, and for now we simply assume that $k_n \xrightarrow{n\rightarrow \infty} \infty$. 

We define a sequence of pairs of finite words $\{\beta_{n,0},\beta_{n,1}\}_{n=1}^{\infty}$ over the alphabet $\Sigma$ of lengths $|\beta_{n,0}|=|\beta_{n,1}|=D_n$. We take $$\beta_{1,0} = 1, \ \beta_{1,1}=2$$
(length-$1$ words) and
\begin{eqnarray*}
    \beta_{n+1,0} & = & \beta_{n,t_{k_n}(0)} \ \beta_{n,t_{k_n}(1)} \ \cdots \ \beta_{n,t_{k_n}(2^{k_n}-1)} \ \underbrace{0 \ \cdots \ 0}_{2^{k_n} D_n\ \text{times}} \\
    \beta_{n+1,1} & = & \tau(\beta_{n+1,0}).
\end{eqnarray*}
Namely, in each step, we use De Bruijn sequences to determine the indices for the sequences from the previous step that will be concatenated to obtain the words in the current one. 

\begin{lem} \label{lem:basic}
\begin{enumerate}
    \item For all $n\geq 1$, we have $|\beta_{n,0}|=|\beta_{n,1}|=D_n$.
    \item For all $n\geq 2$, the words $\beta_{n,0},\beta_{n,1}$ start with $1$ or $2$, and end with a segment of length $\frac{D_n}{2}$ of $0$.
    \item For all $n\geq 1$, the word $\beta_{n,0}$ is a prefix of $\beta_{n+1,0}$ and the word $\beta_{n,1}$ is a prefix of $\beta_{n+1,1}$. In particular, $\beta_{n,0}\neq \beta_{n,1}$.
\end{enumerate}
\end{lem} 

\begin{proof}   
\begin{enumerate}
    \item For $n=1$, this follows since by definition $\beta_{1,0}=1,\ \beta_{1,1}=2$ and $D_1=1$. By induction, $$|\beta_{n+1,1}| = |\beta_{n+1,0}|=2^{k_n} D_n+2^{k_n}D_n=D_{n+1}.$$
    \item By definition, for all $n\geq 2$ we have that $\beta_{n,0},\beta_{n,1}$ end with $\underbrace{0 \ \cdots \ 0}_{2^{k_{n-1}} D_{n-1}\ \text{times}}$ and observe that $\frac{D_n}{2}=2^{k_{n-1}} D_{n-1}$.
    \item By definition, $\beta_{n+1,0}$ starts with $\beta_{n,t_{k_n}(0)}=\beta_{n,0}$ since we assume that $t_{k_n}(0)=0$ for all $n$. The result for $\beta_{n,1},\beta_{n+1,1}$ follows by applying $\tau$.
\end{enumerate}
\end{proof}
Therefore, the sequence $\{\beta_{n,0}\}_{n=1}^{\infty}$ converges to a well-defined (infinite) limit word
$$ w := \lim_{n\rightarrow \infty} \beta_{n,0}.$$

\begin{lem} \label{lem:recursive}
Let $u$ be a factor of $w$ of length $|u|\leq 2^{k_n} D_n$. Then $u$ factors one of
$$\underbrace{0\cdots 0}_{2^{k_n} D_n} \beta_{n+1,*}, \ \beta_{n+1,*} \underbrace{0\cdots 0}_{2^{k_n} D_n}$$ where $*\in \{0,1\}$.
\end{lem}

\begin{proof}
We may assume that $u$ is a subword of $w$ of length $2^{k_n}D_n$ exactly. Let $m\geq n+1$ be the minimum index for which $u$ is a subword of one of 
\begin{equation} \label{eq:forms 3.2}
0^{2^{k_m}D_m} \beta_{m,0},\ 0^{2^{k_m}D_m} \beta_{m,1},\ \beta_{m,0} 0^{2^{k_m}D_m},\ \beta_{m,1} 0^{2^{k_m}D_m}.
\end{equation}
Suppose that $m \geq n+2$. Notice that we can write 
\begin{equation} \label{eq:form m-1 3.2}
\beta_{m-1,0}=B0^{2^{k_{m-2}}D_{m-2}}\ \text{and} \ \beta_{m-1,1}=\tau(B)0^{2^{k_{m-2}}D_{m-2}}
\end{equation}
where $|B|=2^{k_{m-2}} D_{m-2}$.
Then $\beta_{m,0},\beta_{m,1}$ take the form
\begin{equation} \label{equation:3}
B_10^{2^{k_{m-2}}D_{m-2}} B_2 0^{2^{k_{m-2}} D_{m-2}} \cdots B_{2^{k_{m-1}}}0^{2^{k_{m-2}}D_{m-2}} 0^{2^{k_{m-1}}D_{m-1}}
\end{equation}
where, for each $1\leq i\leq 2^{k_{m-1}}$, $B_i\in \{B,\tau(B)\}$.
Since $m\geq n+2$ it follows that $|u|=2^{k_n} D_n \leq 2^{k_{m-2}} D_{m-2}$. Therefore, since $u$ factors one of the forms in (\ref{eq:forms 3.2}) and since it is not longer than the zero blocks between any two $B_i$-blocks in (\ref{equation:3}), it follows that in fact $u$ factors one of $$0^{2^{k_{m-2}}D_{m-2}} B,\ 0^{2^{k_{m-2}}D_{m-2}} \tau(B),\ B 0^{2^{k_{m-2}}D_{m-2}},\ \tau(B) 0^{2^{k_{m-2}}D_{m-2}}.$$
Therefore, by (\ref{eq:form m-1 3.2}), $u$ factors one of 
\begin{align*}
  0^{2^{k_{m-1}}D_{m-1}} & \beta_{m-1,0},  \\ 0^{2^{k_{m-1}}D_{m-1}} & \beta_{m-1,1}, \\ & \beta_{m-1,0} 0^{2^{k_{m-1}}D_{m-1}}, & \\ & \beta_{m-1,1} 0^{2^{k_{m-1}}D_{m-1}},  
\end{align*}
contradicting the minimality of $m$. Hence $u$ factors one of the forms (\ref{eq:forms 3.2}) for $m=n+1$, as claimed.
\end{proof}

\begin{lem} \label{lem:j}
For every $n\gg 1$ and for every $1\leq j\leq k_n$, we have $$p_w(jD_n) \leq (3\cdot 2^{j-1} + 4j) D_n + 4.$$
\end{lem}

\begin{proof}
Let $1\leq j\leq k_n$ be given and consider any factor $u$ of $w$ of length $j D_n$. Since $jD_n \leq 2^{k_n} D_n$, it follows from Lemma \ref{lem:recursive} that $u$ factors one of $0^{2^{k_n} D_n}\beta_{n+1,*}$ or $\beta_{n+1,*}0^{2^{k_n}  D_n}$. Recall that 
$$\beta_{n+1,0}=\beta_{n,t_{k_n}(0)}\beta_{n,t_{k_n}(1)}\cdots \beta_{n,t_{k_n}(2^{k_n}-1)} 0^{2^{k_n} D_n},$$ 
Hence $u$ takes either the form 
\begin{itemize}
    \item $0^i u'$ where $0\leq i\leq jD_n$ and $u'$ is a prefix of $\beta_{n+1,*}$ of length $jD_n - i$;
    \item the form $u'' 0^i$ where $0\leq i\leq jD_n$ and $u''$ is a suffix of $\beta_{n,t_{k_n}(0)}\cdots \beta_{n,t_{k_n}(2^{k_n}-1)}$ or $\beta_{n,\overline{t_{k_n}(0)}}\cdots \beta_{n,\overline{t_{k_n}(2^{k_n}-1)}}$ 
    of length $j D_n - i$; 
    or
    \item $u$ is a factor of $\beta_{n,t_{k_n}(0)}\beta_{n,t_{k_n}(1)}\cdots \beta_{n,t_{k_n}(2^{k_n}-1)}$ or $\beta_{n,\overline{t_{k_n}(0)}}\beta_{n,\overline{t_{k_n}(1)}}\cdots \beta_{n,\overline{t_{k_n}(2^{k_n}-1)}}$. 
\end{itemize}
The first two cases give at most $2(jD_n+1)$ factors each one (corresponding to the value of $i$ and whether $*=0,1$).

In the last case, suppose $u$ factors $$\beta_{n,t_{k_n}(0)}\beta_{n,t_{k_n}(1)}\cdots \beta_{n,t_{k_n}(2^{k_n}-1)},$$ 
whose factors of length $j D_n$ are contained in concatenations of $j+1$ consecutive length-$D_n$ blocks. Therefore, $u$ is contained in a factor of the form $\beta_{n,t_{k_n}(p)}\cdots \beta_{n,t_{k_n}(p+j)}$ for some $0\leq p\leq 2^{k_n}-j-1$. We have at most $2^{j+1}$ such sequences of blocks, and $u$ is determined within such a sequence by its starting position in the first block $\beta_{n,t_{k_n}(p)}$, for which there are $D_n$ possibilities. However, if the starting position of $u$ is within the second half of $\beta_{n,t_{k_n}(p)}$, which is a sequence of zeros regardless of whether $t_{k_n}(p)$ is $0$ or $1$, it does not matter whether the first block is $\beta_{n,0}$ or $\beta_{n,1}$, and hence we need to count only $2^j$ sequences of blocks. Altogether, for the third case, we have at most
\[
\underbrace{\frac{D_n}{2}}_{\substack{\text{Starting position} \\ \text{of}\ u \\ \text{in the first} \\ \text{half of}\ \beta_{n,t_{k_n}(p)}} } \cdot 2^{j+1} + \underbrace{\frac{D_n}{2}}_{\substack{\text{Starting position} \\ \text{of}\ u \\ \text{in the second} \\ \text{half of}\ \beta_{n,t_{k_n}(p)}} } \cdot 2^j
\]
that is, $3\cdot 2^{j-1}D_n$.
The same bound also covers the case that $u$ factors $\beta_{n+1,1}$ rather than $\beta_{n+1,0}$, which we assumed for clarity of the argument.

Combining with the first two cases, we have at most
\[
2(jD_n+1)+2(jD_n+1)+3\cdot 2^{j-1}D_n = (3\cdot 2^{j-1} + 4j) D_n + 4
\] 
possibilities for $u$.
\end{proof}

\begin{lem} \label{lem:difference derivative}
For all $n\gg 1$, we have $$p_w((k_n + 1) D_n) - p_w((k_n - 3) D_n) \geq 2^{k_n-3} D_n.$$
\end{lem}

\begin{proof}
We first argue that
\begin{equation} \label{eq:first}
p_w((k_n + 1) D_n) \geq \frac{(2^{k_n} - k_n) D_n}{2}.
\end{equation}
For each $0\leq r\leq \frac{D_n}{2}$ and $1\leq i\leq 2^{k_n}-k_n$, consider the word
$$\xi_{r,i} := \underbrace{0\ \cdots \ 0}_{r\ \text{times}} \ \beta_{n,t_{k_n}(i)} \ \beta_{n,t_{k_n}(i+1)}\ \cdots \ \beta_{n,t_{k_n}(i+k_n-1)}$$
First, notice that $\beta_{n,t_{k_n}(i-1)}$ takes the form $*\cdots * \underbrace{0\cdots 0}_{D_n/2\ \text{times}}$ by Lemma \ref{lem:basic} (we may assume that $n\geq 2$). Hence $\xi_{r,i}$ factors $$\beta_{n,t_{k_n}(i-1)}\beta_{n,t_{k_n}(i)} \ \beta_{n,t_{k_n}(i+1)}\ \cdots \ \beta_{n,t_{k_n}(i+k_n-1)}$$ which, in turn, factors $\beta_{n+1,0}$ and hence $w$. Next, for each $0\leq r\leq \frac{D_n}{2}$ and $1\leq i\leq 2^{k_n}-k_n$, the length of $\xi_{r,i}$ is $$r+k_n \cdot |\beta_{n,*}|=r+k_n D_n \leq \frac{D_n}{2}+k_n D_n$$ so we can fix an arbitrary right prolongation $\xi_{r,i}u_{r,i}$ of $\xi_{r,i}$, of length $k_n D_n + D_n$. Suppose that $\xi_{r,i} u_{r,i} = \xi_{s,j} u_{s,j}$. Without loss of generality, suppose that $r\leq s$. Then $$|\xi_{r,i}|=r+k_n D_n\leq s+k_n D_n=|\xi_{s,j}|,$$ so $\xi_{r,i}$ is a prefix of $\xi_{s,j}$. Since $\beta_{n,0},\beta_{n,1}$ start with either $1$ or $2$ by Lemma \ref{lem:basic}, the prefix of consecutive zeros in both $\xi_{r,i} u_{r,i}=\xi_{s,j} u_{s,j}$ is equal to both $r$ and to $s$, so $r=s$. Therefore 
$$\beta_{n,t_{k_n}(i)} \ \beta_{n,t_{k_n}(i+1)}\ \cdots \ \beta_{n,t_{k_n}(i+k_n-1)} = \beta_{n,t_{k_n}(j)} \ \beta_{n,t_{k_n}(j+1)}\ \cdots \ \beta_{n,t_{k_n}(j+k_n-1)}$$
and since all $\beta_{n,*}$ are of the same length ($D_n$), it follows that $$t_{k_n}(i) t_{k_n}(i+1) \cdots t_{k_n}(i+k_n-1) = t_{k_n}(j) t_{k_n}(j+1) \cdots t_{k_n}(j+k_n-1)$$ so since $t_{k_n}$ represents a De Bruijn sequence, $i=j$. We thus found at least $\frac{D_n}{2}\cdot (2^{k_n} - k_n)$ distinct factors of $w$ of length $(k_n+1)D_n$, namely, $$\{\xi_{r,i} u_{r,i}\ |\ 0\leq r\leq D_n/2,\ 1\leq i\leq 2^{k_n}-k_n\},$$ and (\ref{eq:first}) follows.

Next, substituting $j=k_n-3$ (here we assume $n\gg 1$ such that $k_n\geq 4$) in Lemma \ref{lem:j}, we get
\begin{equation*}
p_w((k_n - 3) D_n) \leq (3\cdot 2^{k_n-4} + 4 k_n - 12) D_n + 4.
\end{equation*}
If $n\gg 1$ such that $k_n$ is large enough then it follows that
\begin{equation} \label{eq:second}
p_w((k_n - 3) D_n) \leq 2^{k_n-2} D_n.
\end{equation}

Finally, from (\ref{eq:first}) and (\ref{eq:second}), we obtain (as long as $k_n\geq 7$, which we may readily assume)
\begin{eqnarray*}
p_w((k_n+1)D_n) - p_w((k_n-3)D_n) & \geq & \frac{(2^{k_n} - k_n) D_n}{2} - 2^{k_n-2}\cdot D_n \\ & = & (2^{k_n-2} - \frac{1}{2}k_n)D_n \\ & \geq & 2^{k_n - 3} D_n,
\end{eqnarray*}
completing the proof.
\end{proof}

\begin{lem} \label{lem:bounding some p}
We have $p_w(2^{k_n}D_n) \leq 2^{k_n + 4} D_n$.
\end{lem}
\begin{proof}
Let $u$ be a factor of $w$ of length $2^{k_n} D_n$.
By Lemma \ref{lem:recursive}, $u$ factors one of
$$\underbrace{0\cdots 0}_{2^{k_n} D_n} \beta_{n+1,0},\ \ \underbrace{0\cdots 0}_{2^{k_n} D_n} \beta_{n+1,1}, \ \ \beta_{n+1,0} \underbrace{0\cdots 0}_{2^{k_n} D_n}, \ \ \beta_{n+1,1} \underbrace{0\cdots 0}_{2^{k_n} D_n}.$$
Recall that $|\beta_{n+1,*}|=D_{n+1}=2^{k_n+1}D_n$. Hence each one of the above four words is of length $2^{k_n} D_n + 2^{k_n + 1}D_n=3\cdot 2^{k_n} D_n$. Hence
$$p_w(2^{k_n} D_n) \leq 4 \cdot 3\cdot2^{k_n} D_n \leq 2^{k_n+4} D_n,$$
as claimed.
\end{proof}

We now specify our choice of $\{k_n\}_{n=1}^{\infty}$. Let $a>1$ be given; it will be convenient to denote $a=1+\varepsilon$. Recall that $D_1=k_1=1$ and $D_{n+1} = 2^{k_n + 1} D_n$. For every $n\geq 1$, given $D_1,\dots,D_{n+1}$ and $k_1,\dots,k_{n}$, we take $k_{n+1}$ to be the minimum positive integer $x$ such that
$$ D_{n+1}^\varepsilon < \frac{2^{x+1}}{(x+1)^a}$$
(which is well-defined, since the function $f(x)=\frac{2^x}{x^a}$ tends to infinity.) 
Notice that for $n\gg 1$ -- say, $n \geq n_0$ for some $n_0$ depending solely on $\varepsilon$ -- this means that
$$ \frac{2^x}{x^a} \leq D_{n+1}^\varepsilon < \frac{2^{x+1}}{(x+1)^a}$$
We proceed inductively and obtain a sequence $\{k_n\}_{n=1}^{\infty}$ tending to infinity such that for $n\geq n_0$
$$\frac{2^{k_n}}{k_n^a} \leq D_n^\varepsilon < \frac{2^{k_n + 1}}{(k_n + 1)^a}.$$
Moreover, $\{k_n\}_{n=1}^{\infty}$ is non-decreasing, eventually increasing. Once the sequences $\{k_n\}_{n=1}^{\infty}$ and $\{D_n\}_{n=1}^{\infty}$ are determined, the construction proceeds as described above.

\begin{exmpl}[{The construction for $a=2$.}]
Here $n_0=4$. Consider the eventually increasing sequence $\frac{2^t}{t^2}$ for $t=1,2,3,\dots$
$$2,\ 1,\ \frac{8}{9},\ 1,\ \frac{32}{25},\ \frac{64}{36},\ \frac{128}{49},\ \frac{256}{64}, \ \frac{512}{81}, \dots$$

We have $k_1=1$ and, for each $n$, assuming that $k_n$ -- and thus $D_{n+1}=2^{k_n + 1} D_n$ -- have been defined, take $k_{n+1}=t$ the largest possible such that
$$\frac{2^t}{t^2} \leq D_{n+1} < \frac{2^{t+1}}{(t+1)^2}.$$
For instance,
\begin{align*}
    & D_1 = 1, \ \ \ k_1 = 1 \\
    & D_2 = 4, \ \ \ k_2 = 8 \\
    & D_3 = 2^{8+1}\cdot 4=2^{11}, \ \ \ k_3 = 19\\
    & D_4 = 2^{19+1}\cdot 2^{11}=2^{31}, \ \ \ k_4 = 41 \\
    & \vdots
\end{align*}
\end{exmpl}

\begin{figure}[ht]
\centering

\begin{minipage}{0.10\textwidth}
    \centering
    \includegraphics[width=\linewidth]{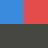}
\end{minipage}
\hspace{0.05\textwidth}
\begin{minipage}{0.10\textwidth}
    \centering
    \includegraphics[width=\linewidth]{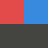}
\end{minipage}
\caption{\footnotesize{A visualization ($a=2$) of $\beta_{2,0}=1200$ (left) and $\beta_{2,1}=2100$ (right). Blue squares represent $1$, red squares represent $2$, and gray squares represent $0$.}}
\end{figure}

\begin{figure}[h]
\centering

\begin{minipage}{0.20\textwidth}
    \centering
    \includegraphics[width=\linewidth]{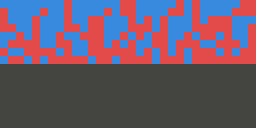}
\end{minipage}
\hspace{0.07\textwidth}
\begin{minipage}{0.20\textwidth}
    \centering
    \includegraphics[width=\linewidth]{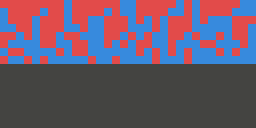}
\end{minipage}
\caption{\footnotesize A visualization ($a=2$) of $\beta_{3,0}$ (left) and $\beta_{3,1}$ (right). Blue squares represent $\beta_{2,0}$, red squares represent $\beta_{2,1}$, and gray squares represent $0$ blocks of length $4$.}
\end{figure}

\begin{figure}[h]
    \centering
    \includegraphics[width=0.64\textwidth]{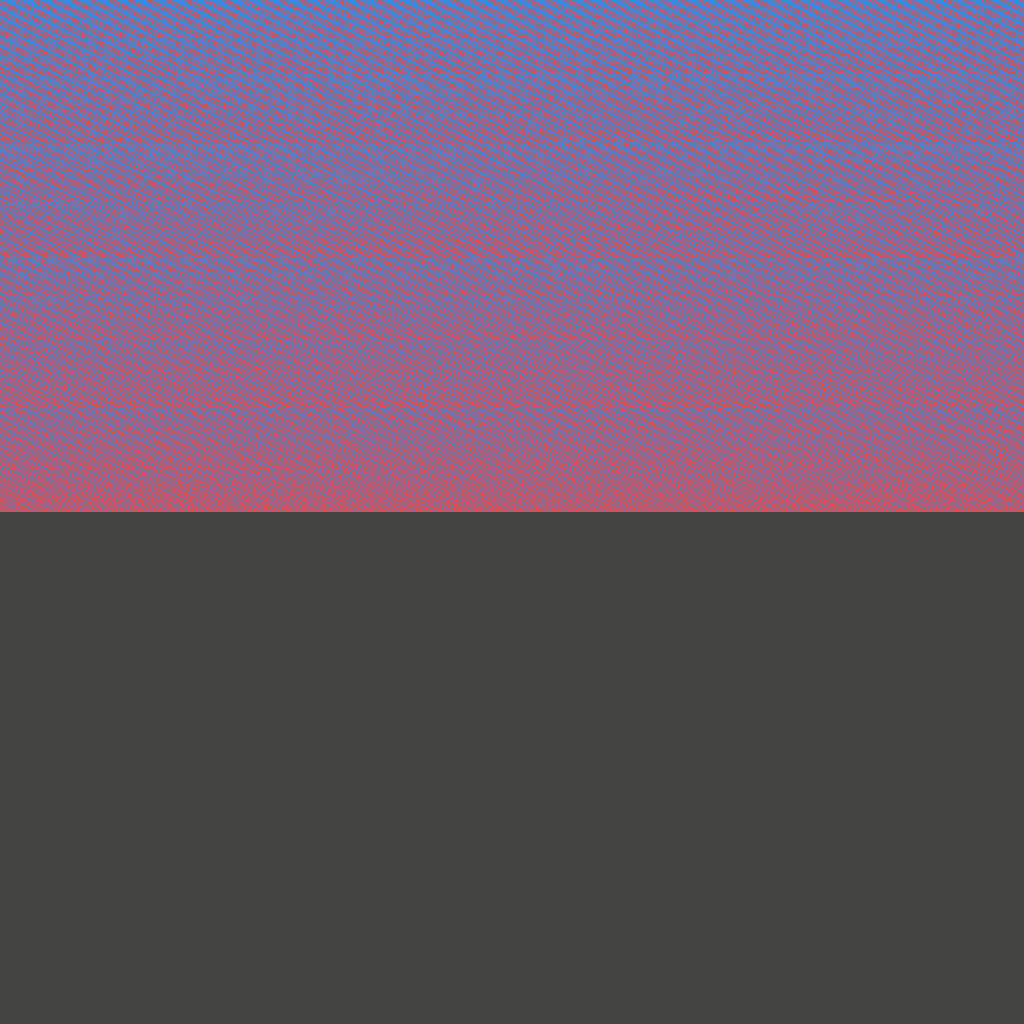}
    \caption{\footnotesize A visualization ($a=2$) of $\beta_{4,0}$. Blue squares represent $\beta_{3,0}$, red squares represent $\beta_{3,1}$, and gray squares represent length-$2048$ blocks of~$0$.}
\end{figure}

\bigskip
\bigskip
\bigskip
\bigskip
\bigskip
\bigskip

Under the above specification,

\begin{prop} \label{prop:quad cplx}
We have $p_w(n)=O(n^a)$.
\end{prop}
\begin{proof}
We assume $n\gg 1$ so that the assertions of the previous lemmata hold, and such that $n \geq n_0$. 

Suppose that $n_0 D_n \leq m \leq k_n D_n$ (we may readily assume that $n\geq n_0$ and that $k_n > n_0$). Let $n_0 \leq j\leq k_n - 1$ be an integer such that $j D_n\leq m \leq (j+1) D_n$. Then, by Lemma \ref{lem:j}, 
\begin{align*}
p_w(m) & \leq p_w((j+1) D_n) \\ & \leq (3\cdot 2^{j}+4j+4) D_n + 4 \\ & \leq  2^{j+3} D_n
\end{align*}
(where the last inequality holds since we may assume $n\gg 1$) and since $n_0 \leq j\leq k_n$, we have 
$$\frac{2^j}{j^a} \leq \frac{2^{k_n}}{k_n^a} \leq D_n^\varepsilon$$
and, together with the previous inequality,
\begin{eqnarray*}
p_w(m) & \leq & 2^{j+3} D_n \\ & \leq & 8 \cdot (j^a D_n^\varepsilon) D_n \\ & = & 8 \cdot (j D_n)^a \\ & \leq & 8 m^a.
\end{eqnarray*}

Now suppose that $k_n D_n \leq m \leq 2^{k_n} D_n$. Since $\frac{2^{k_n}}{k_n^a} \leq D_n^\varepsilon$, we have $2^{k_n} D_n\leq k_n^a D_n^a$ and therefore, by Lemma \ref{lem:bounding some p},
\begin{eqnarray*}
p_w(m) & \leq & p_w(2^{k_n} D_n) \\ & \leq & 2^{k_n+4} D_n \\ & \leq & 16\cdot k_n^a D_n^a\leq 16 m^a.
\end{eqnarray*}

Finally, suppose that $2^{k_n} D_n \leq m \leq n_0 D_{n+1}$. 
Then $$2n_0 m \geq 2n_0 \cdot 2^{k_n} D_n = n_0 D_{n+1}$$ and $$2n_0 m \leq 2n_0^2 D_{n+1} \leq k_{n+1} D_{n+1}$$ where the last inequality holds if $n\gg 1$ is large enough to ensure that $k_{n+1} > 2n_0^2$. Now by the first case treated in this proof,
$$p_w(m) \leq p_w(2n_0 m) \leq 8 \cdot (2n_0 \cdot m)^a = 8\cdot 2^a n_0^a \cdot m^a.$$
Hence, for all $m\gg 1$ we have that $$p_w(m) \leq K \cdot m^a,$$ for $K=16 \cdot 2^a\cdot n_0^a$, so $p_w(m)=O(m^a)$.
\end{proof}

\begin{proof}[{Proof of Theorem \ref{thm:quadratic}}]
By Proposition \ref{prop:quad cplx}, $p_w(n)=O(n^{1+\varepsilon})$. By Lemma \ref{lem:difference derivative}, for every $n\gg 1$, we have that
$$p_w((k_n+1)D_n) - p_w((k_n-3)D_n) \geq 2^{k_n - 3}D_n.$$
Since 
$$p_w((k_n+1)D_n) - p_w((k_n-3)D_n) = \sum_{r=(k_n-3)D_n + 1}^{(k_n+1)D_n} p_w'(r),$$
we conclude that for some $r \in [(k_n-3)D_n + 1,(k_n+1)D_n]$,
\begin{align*}
p_w'(r) & \geq \frac{1}{4D_n} \sum_{s=(k_n-3)D_n+1}^{(k_n+1)D_n} p_w'(s) \\ & = \frac{p_w((k_n+1)D_n) - p_w((k_n-3)D_n)}{4D_n} \\ & \geq 2^{k_n - 5}.    
\end{align*}
Recall that $$D_n^\varepsilon < \frac{2^{k_n+1}}{(k_n+1)^a},$$ so $2^{k_n} > \frac{1}{2} D_n^\varepsilon (k_n+1)^a$ and hence
\begin{align} \label{eq:final}
p_w'(r) & \geq 2^{k_n-5} \\ & \geq 2^{-6}\cdot (k_n+1)^a D_n^\varepsilon \nonumber \\ & \geq \frac{k_n + 1}{64}\cdot r^\varepsilon. \nonumber
\end{align}
Finally,
\begin{align*} 
r & \leq (k_n+1) D_n \\ & \leq (k_n+1)\cdot \frac{2^{(k_n+1)/\varepsilon}}{(k_n + 1)^{a/\varepsilon}} \\ & = \frac{2^{(k_n+1)/\varepsilon}}{(k_n + 1)^{1/\varepsilon}} \\ & \leq (2^{1/\varepsilon})^{k_n+1}
\end{align*}
so $k_n+1 \geq \varepsilon\cdot \log_2 r$ and, together with (\ref{eq:final}), $$p_w'(r) \geq \frac{\varepsilon}{64} \cdot r^\varepsilon \log_2 r$$ for infinitely many values of $r$, and therefore $p_w'(n) \neq O(n^\varepsilon)$.
\end{proof}

\section{Polynomial complexity and exponential-type derivatives}

\subsection{Infinite words}

We work over the alphabet $$\Sigma = \{0,1,2\}.$$ 
We define sequences $\{n_k\}_{k=1}^{\infty}$ and $\{s_k\}_{k=1}^{\infty}$ by 
$$n_1=1,\ n_2=2,\ n_3 = 4,\ n_4 = 8,\ n_5=16, \ \text{and}\ \ n_k = \lceil \frac{2^{2^{k-2}}}{k^2} \rceil \ \ \text{for}\ \ k\geq 6$$
and
$$s_k = 2^{2^{k-1}} \ \ \text{for}\ \ k\geq 1.$$
Notice that $n_{k+1}\geq 2n_k$ for all $k\geq 1$, and $n_{k+1}\geq 3n_k$ for all $k\geq 5$.
Define a sequence of words over $\Sigma$ as follows. We take $$ X_1 = \{1,2\}$$ and, for each $k\geq 1$, $$X_{k+1} = X_k \ 0^{n_{k+1}-2n_k} \ X_k = \{u0^{n_{k+1}-2n_k} v\ |\ u,v\in X_k\}.$$
For instance,
\begin{align*}
    & X_1=\{1,2\},\ X_2=\{1,2\}^2,\ X_3=\{1,2\}^4,\ X_4=\{1,2\}^8,\ X_5=\{1,2\}^{16},\\ & X_6 = \{1,2\}^{16} \underbrace{0\dots 0}_{1789\ \text{times}} \{1,2\}^{16}, \\
    & X_7 = \{1,2\}^{16} \underbrace{0\dots 0}_{1789\ \text{times}} \{1,2\}^{16} \underbrace{0\dots 0}_{\substack{87648752 \\ \text{times}}} \{1,2\}^{16} \underbrace{0\dots 0}_{1789\ \text{times}} \{1,2\}^{16}.
\end{align*}
Notice that for every $k\geq 1$ we have $|X_k| = s_k$; indeed, by induction, $$|X_{k+1}| = |X_k|^2=s_k^2=s_{k+1}$$ and $X_k\subseteq \Sigma^{n_k}$ by construction. Notice that for every $1\leq j\leq k$, we can write $X_k = S\cdot X_j$ for some set $S$. Hence if $i\leq n_k-n_j$, we have that $\text{pre}^{n_k}_i(X_k)$ is an image (under a prefix map) of the set $S$, so 
\begin{equation*}
|\text{pre}^{n_k}_i(X_k)|\leq |S|=\frac{|X_k|}{|X_j|}=\frac{s_k}{s_j}.
\end{equation*}

Moreover, for every $\alpha,\beta \in \bigcup_{k=1}^{\infty} X_k$ we have a word of the form $\alpha \cdots \beta$ factoring a word from some $X_l$; hence the hereditary language $\mathcal{L}=\bigcup_{n=0}^{\infty} \mathcal{L}_n$ of $\bigcup_{k=1}^{\infty} X_k$ coincides with the hereditary language of some recurrent infinite word $w\in \Sigma^\mathbb{N}$, and thus, with the hereditary language associated with a recurrent bi-infinite word.

\begin{lem} \label{lem:1.1 complexity above}
    There exists some constant $K>0$ such that  $p_w(n_k)\leq K \cdot 2^{2^{k-1}}$ for all $k\geq 1$.
\end{lem}
\begin{proof}
    Evidently, it suffices to prove the claim for all $k\geq 5$.
    We first claim that every word $u\in \mathcal{L}_{n_k},\ k\geq 5$ factors some word from $0^{n_k}X_k\cup X_k 0^{n_k}$. Suppose that $m$ is the minimum index for which $u$ factors a word from $0^{n_m}X_m\cup X_m 0^{n_m}$ and suppose that $m>k$. Since $$X_m=X_{m-1}0^{n_m-2n_{m-1}} X_{m-1},$$ it follows that $u$ factors a word from
    \begin{equation*} 
    0^{n_m} X_{m-1}0^{n_m-2n_{m-1}} X_{m-1} \cup X_{m-1}0^{n_m-2n_{m-1}} X_{m-1} 0^{n_m}.
    \end{equation*}
    Since $X_m\subseteq \Sigma^{n_m}$ and $n_m-2n_{m-1}\geq n_{m-1}\geq n_k$ (as $m\geq 6$) we conclude that $u$ factors a word from $$0^{n_{m-1}} X_{m-1} \cup X_{m-1} 0^{n_{m-1}},$$ contradicting the minimality of $m$. Hence $u$ factors a word from $0^{n_k} X_{k} \cup X_{k} 0^{n_k}$.

    Let us analyze the number of length-$n_k$ subwords of $0^{n_k} X_k$. If $u$ is such a word then $u=0^{n_k-i}v$ for $v$ an $i$-prefix of a word from $X_k$, namely, $v\in \text{pre}^{n_k}_i(X_k)$ for some $0\leq i\leq n_k$. Suppose $n_k - n_{j+1} < i \leq n_k - n_j$. Since $X_k$ can be written as $S\cdot X_j$, we have at most $\frac{|X_k|}{|X_j|}=\frac{s_k}{s_j}$ options for such $u$.

    \begin{align*}
    \Big|\Big\{\substack{\text{Subwords of}\ 0^{n_k} X_k \\ \text{of length}\ n_k}\Big\}\Big| & \leq \sum_{i=0}^{n_k}|\{0^{n_k-i}v|v\in \text{pre}^{n_k}_i(X_k)\}| \\ & = \sum_{i=0}^{n_k} |\text{pre}^{n_k}_i(X_k)| \\ & \leq \underbrace{|X_k|}_{i=n_k} + \underbrace{1}_{i=0} + \sum_{j=1}^{k-1} \underbrace{(n_{j+1} - n_j)}_{n_k-n_{j+1}<i\leq n_k-n_j} \cdot \frac{s_k}{s_j} \\  & = 1 + s_k \cdot \left(1 + \sum_{j=1}^{k-1} \frac{n_{j+1} - n_j}{s_j} \right) \\ & \leq 1 + s_k \cdot \left( 1 + \sum_{j=1}^{\infty} \frac{n_{j+1}}{s_j} \right) \\ & \leq 1 + s_k \cdot \left( 1 + n_2 + n_3 + n_4 + n_5 + \sum_{j=5}^{\infty} \frac{n_{j+1}}{s_j} \right) \\ & \leq 1 + s_k \cdot \left(31 + \sum_{j=1}^{\infty} \frac{\frac{2^{2^{j-1}}}{j^2}+1}{2^{2^{j-1}}} \right) \\ & \leq 1 + s_k \cdot \left( 31 + \sum_{j=1}^{\infty} \frac{1}{j^2} + \sum_{j=1}^{\infty} 2^{-2^{j-1}} \right)\leq C\cdot s_k
    \end{align*}
    for some constant $C>0$, as the series in the penultimate expression are all convergent.
    
    The analysis for subwords of $X_k 0^{n_k}$ of length $n_k$ is symmetric, so altogether the number of words of length $n_k$ factoring $0^{n_k} X_k \cup X_k 0^{n_k}$ is at most $2C\cdot s_k$, so $p_w(n_k) \leq 2C\cdot s_k = K \cdot 2^{2^{k-1}}$.
\end{proof}

\begin{prop} \label{prop:cplx 1}
    We have $p_w(n)\leq n^{3+o(1)}$.
\end{prop}
\begin{proof}
    Let $k\geq 6$ and suppose that $n_k\leq t \leq n_{k+1} - 2n_k$. We claim that every word in $\mathcal{L}_t$ takes one of the forms 
    \begin{align} \label{eq:forms in prop 1}
     0^i v 0^{t-i-n_k} & \ \ \ \text{for}\ 0\leq i\leq t-n_k,\ v\in X_k \\
     0^i \text{pre}^{n_k}_{t-i}(v) & \ \ \ \text{for}\ t-n_k\leq i\leq t,\ v\in X_k \nonumber \\
     \text{suf}^{n_k}_{t-i}(v) 0^i & \ \ \ \text{for}\ t-n_k\leq i\leq t,\ v\in X_k. \nonumber
    \end{align}
    As argued in the proof of Lemma \ref{lem:1.1 complexity above}, every word in $\mathcal{L}$ of length $n_{k+1}$ factors a word from $0^{n_{k+1}} X_{k+1} \cup X_{k+1} 0^{n_{k+1}}$ and thus factors a word from 
    \begin{equation} \label{eq:set prop 1}    0^{n_{k+1}}\underbrace{X_k0^{n_{k+1}-2n_k}X_k}_{X_{k+1}} \cup \underbrace{X_k0^{n_{k+1}-2n_k}X_k}_{X_{k+1}}0^{n_{k+1}}.
    \end{equation}
    Pick a word $u\in \mathcal{L}_{t}$. It factors some word in $\mathcal{L}_{n_{k+1}}$ and therefore factors a word from (\ref{eq:set prop 1}). Moreover, since $u$ is of length $t \leq n_{k+1}-2n_k$, it follows that $u$ factors a word from $0^\infty X_k 0^\infty$. It follows that $u$ takes one of the forms listed in (\ref{eq:forms in prop 1}). The number of words of the form $0^i v 0^{t-i-n_k}$ is at most $(t-n_k+1)|X_k|$; the number of words of the form $0^i\text{pre}^{n_k}_{t-i}(v)$ is at most $(t+1)|X_k|$; and likewise the third form listed in (\ref{eq:forms in prop 1}). Altogether, we have at most $$3(t+1)|X_k|=3(t+1)s_k.$$ Notice that for all $k\geq 5$, $$n_k\geq \frac{2^{2^{k-2}}}{k^2} = \frac{\sqrt{s_k}}{k^2}$$ so $s_k\leq k^4 n_k^2$. Moreover, $$n_k \geq \frac{2^{2^{k-2}}}{k^2}\geq 2^{2^{k-3}}$$ for all $k\geq 6$, so $k\leq \log \log n_k + 3$.   
    Hence \begin{align*}|\mathcal{L}_t| & \leq 3(t+1)\cdot k^4 n_k^2 \\  & \leq 3(t+1)n_k^{2}(\log \log n_k +3)^4 \\ & \leq 3(t+1)t^2 (\log \log t + 3)^4.\end{align*} This is valid for any $n_k\leq t \leq n_{k+1}-2n_k$. Since $n_{k+1}-2n_k \geq \frac{1}{2}n_{k+1}$, for every $n\gg 1$, either $n$ or $2n$ belong to an interval of the form $[n_k,n_{k+1}-2n_k]$. It follows that $p_w(n)\leq n^{3+o(1)}$.
\end{proof}

\begin{lem} \label{lem:1.1 derivative below}
    We have $p_w'(n_k) \geq 2^{2^{k-1}-1}$ for all $k\geq 2$.
\end{lem}
\begin{proof}
By induction, for every $k\geq 1$, we have that $X_{k+1}=S_{k+1}\cdot \{1,2\}$ and $$|S_{k+1}|=\frac{|X_{k+1}|}{2}=\frac{2^{2^{k}}}{2}=2^{2^{k}-1}.$$
Indeed, $X_2=\{1,2\}^2$; and for the induction step, $$X_{k+1}=X_k 0^{n_{k+1}-2n_k}X_k=\underbrace{X_k 0^{n_{k+1}-2n_k} S_k}_{=:S_{k+1}} \{1,2\}.$$
Hence, for every $k\geq 2$, each one of the $|S_k|=2^{2^{k-1}-1}$ words in $S_k$ can be extended to the right by both $1$ and $2$. 
By Lemma \ref{lem:derivative}, $p_w'(n_k)\geq 2^{2^{k-1}-1}$ for all $k\geq 2$.
\end{proof}

\begin{cor} \label{cor:1}
    We have $p_w'(n_k)\geq \frac{1}{K} p_w(n_k)$ for some constant $K>0$ for all $k\geq 2$.
\end{cor}
\begin{proof}
    We have $$p_w'(n_k) \overset{\substack{\text{Lemma} \\ \text{\ref{lem:1.1 derivative below}}}}{\geq} 2^{2^{k-1}-1}\overset{\substack{\text{Lemma} \\ \text{\ref{lem:1.1 complexity above}}}}{\geq} \frac{1}{K} p_w(n_k),$$ as claimed.
\end{proof}

\begin{proof}[{Proof of Theorem \ref{thm:p'>cp}}]

This follows from Proposition \ref{prop:cplx 1} and Corollary \ref{cor:1}.
\end{proof}

\subsection{Subshifts of quasi-linear complexity}
We work over the binary alphabet $\Sigma=\{0,1\}$. Consider the increasing sequence $r_n = \lfloor \frac{2^n}{n^2} \rfloor$ for $n\geq 6$ (i.e. $r_6=1,r_7=2,r_8=4,\dots$). Let $\mathcal{X}$ be the subshift generated by
$$\mathcal{A} = \Big\{\underline{a}=(a_k)_{k=0}^{\infty}\in \{0,1\}^{\mathbb{N}} \ |\ a_k=0\ \forall k\notin \{r_6,r_7,\dots\}\Big\} \subseteq \{0,1\}^{\mathbb{N}}$$
namely, $\mathcal{X}=\overline{\bigcup_{i=0}^{\infty} T^i(\mathcal{A})}$, where $T\colon \{0,1\}^\mathbb{N}\rightarrow \{0,1\}^\mathbb{N}$ is the shift operator.

\begin{lem} \label{lem:der subshift example}
For every $n \geq 7$, we have $p_\mathcal{X}'(r_n)\geq 2^{n-6}$.
\end{lem}

\begin{proof}
Let $n\geq 7$. Consider all words of the form $u=u_1 u_2 \cdots u_{r_n-1}\in \{0,1\}^{r_n-1}$ with $u_k = 0$ for all $k \notin \{r_6,r_7,\dots\}$.
These take the form
$$b_6b_70b_80^{r_9-r_8-1} b_9 \dots b_{n-1} 0^{r_n-r_{n-1}-1}$$
with arbitrary $b_6,\dots,b_{n-1}\in \{0,1\}$, and $u_{r_i}=b_i$. Since the next spot is $r_n$, each of the above words can be prolonged to the right by both $0$ and $1$. Evidently, there are $2^{n-6}$ such words corresponding to all choices of $b_6,\dots,b_{n-1}$. By Lemma \ref{lem:derivative}, we conclude that $p_\mathcal{X}'(r_n)\geq 2^{n-6}$.
\end{proof}

\begin{lem} \label{lem:subshift complexity bound}
For every constant $C\in \mathbb{N}$ there exists some $K=K_C$ such that $p_\mathcal{X}(Cr_n) \leq K \cdot 2^n$ for every $n$.
\end{lem}

\begin{proof}
Fix $C\in \mathbb{N}$; we count factors $u$ of $\mathcal{X}$ of length $Cr_n$. Every such $u$ is a subword of some $\underline{a}\in \mathcal{A}$. 
There exists some $t$ depending only on $C$ such that for all $n$, we have $$Cr_n < r_{i+1}-r_{i},\ \text{for all}\ i\geq n+t.$$ 
Indeed, the difference sequence $\Delta_i := r_i - r_{i-1}$ grows exponentially in $i$ and furthermore, there exists some (universal) $\theta>0$ such that $\Delta_{i+1}\geq \theta r_i$; hence if $i\geq n+t$ then $\Delta_{i+1} \geq \theta r_{n+t}\geq \theta(1+\theta)^t r_n$. If $t$ is large enough depending on $C$ then $\theta(1+\theta)^t > C$.

Therefore, if $u$ is a subword of some $\underline{a}\in \mathcal{A}$ and $u=a_ja_{j+1}\cdots a_{j+Cr_n-1}$ and $j>r_{n+t}$ then $u$ has at most one non-zero letter; there are $1+Cr_{n}$ such words of length $Cr_n$. 

Now fixing $0\leq j\leq r_{n+t}$, the number of possible subwords of some $\underline{a}\in \mathcal{A}$ of the form $a_ja_{j+1}\cdots a_{j+Cr_n-1}$ is $$2^{|\{r_6,r_7,\dots\}\cap [j,j+Cr_n-1]|},$$ i.e., $2$ to the number of `free bits' $b_6,b_7,\dots$ falling in the right window.

Furthermore, $j+Cr_n-1\leq r_{n+t}+Cr_n-1 < r_{n+t+1}$ and therefore
$$|\{r_6,r_7,\dots\}\cap [0,j+Cr_n-1]| \leq n+t-5$$ and
$$|\{r_6,r_7,\dots\}\cap [j,j+Cr_n-1]| \leq n+t-5 - |\{r_6,r_7,\dots\}\cap [0,j-1]|.$$
Hence
\begin{align*}
p_\mathcal{X}(Cr_n) & \leq 1 + Cr_n + \sum_{j=0}^{r_{n+t}} 2^{|\{r_6,r_7,\dots\}\cap [j,j+Cr_n-1]|} \\
& \leq 1 + Cr_n + \sum_{j=0}^{r_{n+t}} 2^{n+t-5-|\{r_6,r_7,\dots\}\cap [0,j-1]|} \\ 
& = 1 + Cr_n + \sum_{k=0}^{n+t-5} 2^{n+t-5-k} \cdot \Big|\{0\leq j\leq r_{n+t}\ : \ |\{r_6,r_7,\dots\}\cap [0,j-1]|=k\}\Big|
\end{align*}
Notice that $|\{r_6,r_7,\dots\}\cap [0,j-1]|=k$ if and only if $r_{k+5} < j \leq r_{k+6}$. It follows that
\begin{align*}
p_\mathcal{X}(Cr_n) & \leq 1 + Cr_n + \sum_{k=0}^{n+t-5} 2^{n+t-5-k} \cdot |\{j\ :\ r_{k+5} < j \leq r_{k+6}\}| \\
& = 1 + Cr_n + 2^{n+t-5} \cdot (r_6 + 1) + \sum_{k=1}^{n+t-5} 2^{n+t-5-k} \cdot (r_{k+6} - r_{k+5}) \\
& =  1 + Cr_n + 2^{n+t-4} + 2^{n+t-5} \cdot \sum_{k=1}^{n+t-5} 2^{-k} \cdot \left(\Big\lfloor \frac{2^{k+6}}{(k+6)^2} \Big\rfloor - \Big\lfloor \frac{2^{k+5}}{(k+5)^2} \Big\rfloor \right) \\
& \leq  1+ Cr_n +  2^{n+t-4} + 2^{n+t-5} \cdot \sum_{k=1}^{n+t-5} 2^{-k} \cdot \frac{2^{k+6}}{(k+6)^2} \\ & \leq 1 + Cr_n +  2^{n+t-4} + 2^{n+t+1}\cdot \sum_{k=1}^{\infty} \frac{1}{k^2} \\ & \leq 1 + \left(\frac{C}{n^2}+2^t\left(\frac{1}{16} + \frac{\pi^2}{3}\right)\right)2^n < K \cdot 2^n
\end{align*}
if we take $K\gg 1$ depending only on $C$.
\end{proof}

\begin{figure}[h]
    \centering \includegraphics[width=0.95\textwidth]{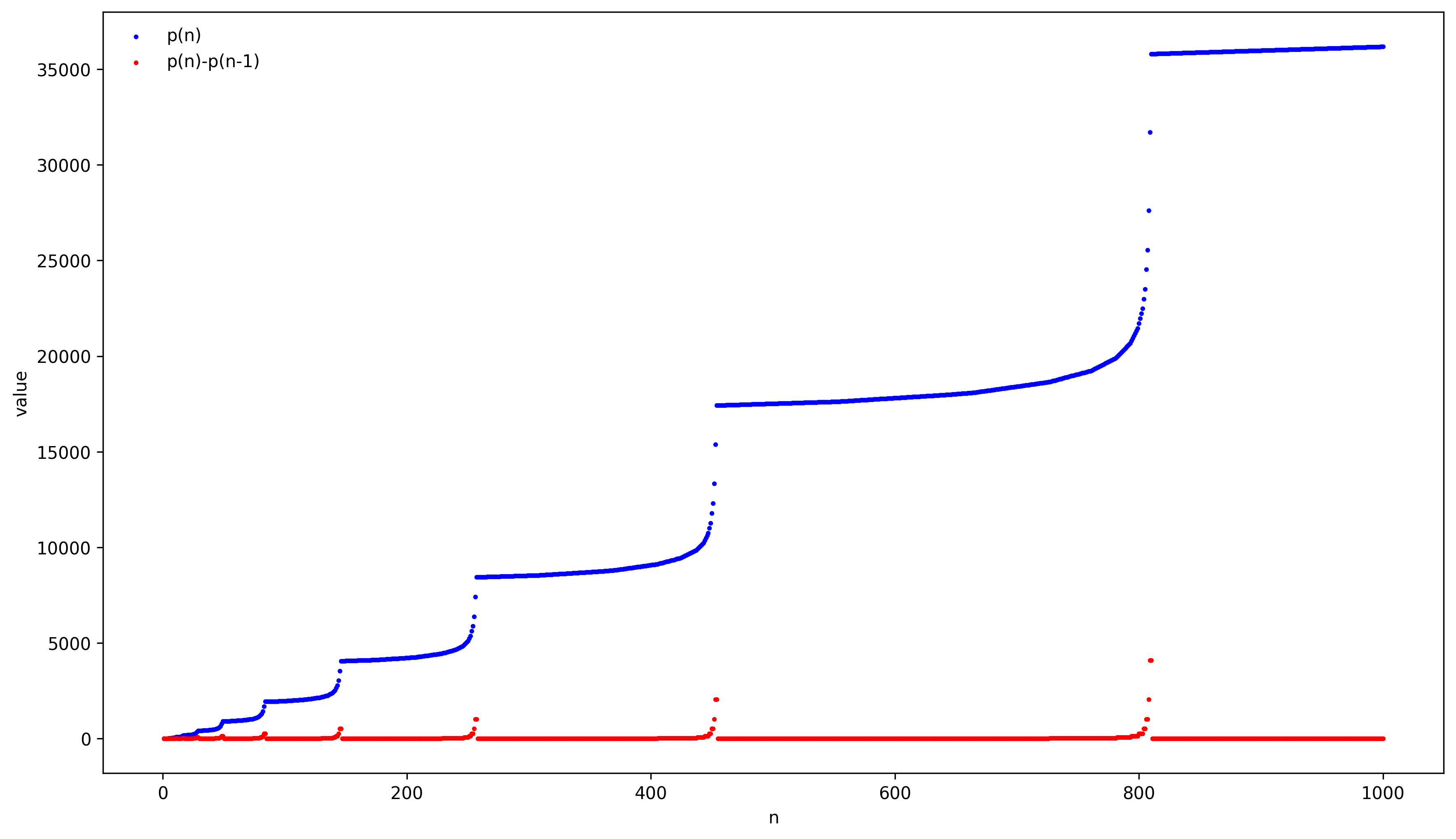}
    \caption{\footnotesize The complexity function $p_\mathcal{X}(n)$ (blue) and its discrete derivative $p_\mathcal{X}'(n)$ (red), for $n\in [1,1000]$.}
\end{figure}

\begin{proof}[{Proof of Proposition \ref{prop:p'>pc}}]
Consider the subshift $\mathcal{X}$ constructed above.
Fix $m$ large enough. Let $n$ be such that $r_n \leq m < r_{n+1}$; we may assume $n\geq 7$, as required by Lemma \ref{lem:der subshift example}. We have $\frac{2^n}{n^2} \leq r_n+1\leq m + 1$ so $2^n\leq n^2(m+1)$, and $n\leq 2\log n + \log (m+1)$, so $n\leq 2\log (m+1)$.
By Lemma \ref{lem:subshift complexity bound}, we have 
\begin{align*}
p_\mathcal{X}(m) \leq p_\mathcal{X}(r_{n+1}) & \leq K\cdot 2^{n+1} \\ & \leq 2K  n^2(m+1) \leq 8K (m+1)\log^2(m+1).
\end{align*} 
By Lemma \ref{lem:der subshift example}, for every $n\geq 7$ we have that $p_\mathcal{X}'(r_n)\geq 2^{n-6}$. By Lemma \ref{lem:subshift complexity bound}, for every $C\in \mathbb{N}$ there exists $K$ such that $p_\mathcal{X}(C r_n) \leq K\cdot 2^n$. It follows that $$p_\mathcal{X}(C r_n) \leq 64K \cdot p_\mathcal{X}'(r_n),$$
as claimed.
\end{proof}

\section{Minimal subshifts of prescribed complexity}

Our goal in this section is to prove Theorem \ref{thm:ergodic}. 

\subsection{Fractal systems of words} \label{subsec:construction}
We start by recalling the construction from \cite{GMZ}. The idea behind it is rooted in the construction of Smoktunowicz and Bartholdi in \cite{SB}, originally in the context of hereditary languages (utilized there for monomial algebras). Building on this self-similar construction, it has been shown in \cite{GaApp,BGJalg,BGIsr,GMZ} how to construct infinite words with strong recurrence properties of well-approximated growth. The state of the art of this evolving line of constructions was given in \cite{GMZ2}, where the first, third and fourth named authors proved that for any non-decreasing, submultiplicative function $f\colon \mathbb{N}\rightarrow \mathbb{N}$ there exists a strictly ergodic (that is, minimal and uniquely ergodic) subshift $\mathcal{X}$ whose complexity function satisfies $f(n) \preceq p_\mathcal{X} \preceq nf(n)$.

\medskip

Let $\Sigma$ be a set of cardinality $b:=f(1)$. We define a sequence of positive integers $\{c_{k}\ \colon\ k=0,1,2,\dots\}$ and define $N_k:=bc_0c_1\cdots c_k$ accordingly.

\begin{prop}[{\cite{GMZ2}}] \label{prop:define}
Given a non-decreasing, submultiplicative, subexponential function $f\colon \mathbb{N}\rightarrow \mathbb{N}$, we can define $\{c_k\ :\ k=0,1,2,\dots\}$ as above such that:
\begin{enumerate}
\item $1\leq c_{k+1}\leq N_k$ for all $k\geq 0$ and $c_0\leq b$
\item $f(2^k)\leq N_k \leq 2f(2^{k+1})$ for all $k\geq 0$
\item $c_k=1$ for infinitely many values of $k$.
\end{enumerate}
\end{prop}

Given a non-decreasing, submultiplicative, subexponential function $f\colon \mathbb{N}\rightarrow \mathbb{N}$, fix a sequence $\{c_k\ \colon\ k=0,1,2,\dots\}$ as in Proposition \ref{prop:define}. Let $\mathcal{S}\subseteq \mathbb{N}$ be the set of indices $k$ for which $c_k=1$, which, as guaranteed, is infinite. 
We define, for each $k=0,1,2,\dots$, sets $$W(k)\subseteq \Sigma^{2^k}, \ C(k)\subseteq W(k),\ U(k)$$ such that $$|C(k)|=c_k, \ |W(k)|=bc_0\cdots c_{k-1}, \  U(k)\subseteq \bigcup_{i=0}^{k} W(i).$$ Moreover, each $U(k)$ is in fact an \emph{ordered list} which can be written as $U(k)=(u^{(k)}_{1},\dots,u^{(k)}_{t_k})$.

We start with $W(0)=\Sigma$ (notice that $|W(0)|=|\Sigma|=b$) and $U(0)$ is a list containing the elements of $W(0)$, arbitrarily ordered. Suppose that $$\{W(i),U(i)\ \colon\ i=0,1,\dots,k\}, \ \{C(i)\ \colon\ i=0,1,\dots,k-1\}$$ have been defined for some $k\geq 0$. 

If $k\notin \mathcal{S}$, or if $k\in \mathcal{S}$ but $2^k<|u^{(k)}_1|$ then we take $C(k)$ to be an arbitrary subset of $W(k)$ of size $c_k$, which is indeed possible since by Proposition \ref{prop:define}, $1 \leq c_k \leq bc_0\cdots c_{k-1}$. We let $W(k+1)=W(k)C(k)$ and update $U(k+1)$ by appending $W(k+1)$, arbitrarily ordered, to the end of $U(k)$.

If $k\in \mathcal{S}$ (so $c_k=1$) and $|u^{(k)}_1| \leq 2^k$ then we fix an arbitrary word $v\in W(k)$ whose prefix is $u^{(k)}_1$, and take $C(k)=\{v\}$. We take $$W(k+1)=W(k)C(k)=W(k)\cdot \{v\}$$ and  update $U(k+1)$ by appending $W(k+1)$, arbitrarily ordered, to the end of $U(k)$, and omitting the first word in it (so $u^{(k+1)}_1=u^{(k)}_2$).

\medskip

Pick an arbitrary infinite word
$$ w \in W(0)C(0)C(1)C(2)C(3)\cdots$$
In \cite{GMZ2} it is proved that $w$ is uniformly recurrent and has well-defined factor frequencies (equivalently, the closure of its shift orbit is strictly ergodic). The hereditary language of finite subwords of $w$ coincides with the hereditary closure of $\bigcup_{k=0}^{\infty} W(k)$.

\subsection{Interleaving words}
Write
$$ w = w_1 w_2 w_3 \dots $$
so $w_{d 2^k+1}\cdots w_{(d+1) 2^k} \in W(k)$ for every $d,k\geq 0$. Moreover, recall that by the above construction (Subsection \ref{subsec:construction}), for every $v\in W(k)$ we have that, for some $l>k$, $C(l)=\{z\}$ where $v$ is the length-$2^k$ prefix of $z$. It follows that 

\begin{cor} \label{cor:every power of 2}
For every $v\in W(k)$ there exists some $l>k$ such that for every $d\geq 1$ odd, we have $$ w_{d 2^l + 1} \cdots w_{d 2^l + 2^k} = v. $$
\end{cor}

\begin{lem} \label{lem:starting position modulo}
Fix $n$ and let $i \in \{1,\dots,2^n\}$ be given. Then there are at most $4f(2^{n+1})$ distinct factors of $w$ of length $2^n$ whose starting position is congruent to $i$ modulo $2^n$; namely, of the form $w_{d 2^n+i}\cdots w_{(d+1) 2^n + i - 1}$.
\end{lem}

\begin{proof}
Observing that $$ w\in W(0)C(0)C(1)C(2)\cdots \subseteq W(n)C(n)W(n)C(n)\cdots $$
we deduce that, for every $d$, $$w_{d 2^n + i} \cdots w_{(d+1) 2^n + i -1} = uv$$ where $|u|=2^n - i + 1,\ |v| = i-1$ and either (i) $u$ is a suffix of a word from $W(n)$ and $v$ is a prefix of a word from $C(n)$, or (ii) $u$ is a suffix of a word from $C(n)$ and $v$ is a prefix of a word from $W(n)$. In either case (recall that $i$ is fixed), $u,v$ have at most $$|W(n)| \cdot |C(n)| = |W(n+1)| = N_n \leq 2f(2^{n+1})$$ options (utilizing Proposition \ref{prop:define}); altogether, we have at most $4f(2^{n+1})$ options for $(u,v)$ and, subsequently, to $w_{d2^n+i}\cdots w_{(d+1)2^n+i-1}$.
\end{proof}

Recall the \emph{Thue-Morse sequence} $$t=t_0 t_1 t_2 \dots = 01101001\dots$$ 
This sequence has many equivalent definitions. For instance, it can be defined as follows: let $t[0]=0$ and $t[n+1]=t[n]\overline{t[n]}$ (where $\overline{u}$ is the bit-wise `NOT' of $u$). Then $t=\lim_{n\rightarrow \infty} t[n]$. Notice that $t[n]=t_0\cdots t_{2^n-1}$ is the $2^n$-prefix of $t$. 
For each $n$, the Thue-Morse sequence takes the form
$$t[n] \ \overline{t[n]}\ \overline{t[n]} \ t[n] \ \cdots$$
namely, the substitution of $t$ under $0\mapsto t[n],1\mapsto \overline{t[n]}$. In particular, $t[n]$ occurs in $t$ with starting position $2^n \cdot j$ for all $j$ such that $t_j=0$ (namely, for all $j$ with an even number of $1$'s in its binary expansion).

An equivalent definition is by parity of the binary expansion: $t_n$ is equal to the number of $1$'s in the binary expression of $n$, modulo $2$.

\begin{lem} \label{lem:syndeticTM}
Fix $n\geq 0$ and $s\geq 0$. Then the set $$\{i\in \mathbb{Z}_{\geq 0}\ |\ t[i,i+2^n-1]=t[n]\ \text{and}\ i=2^{n+s} \cdot \text{odd}\}$$ is \emph{syndetic}, that is, has bounded gaps.
\end{lem}
\begin{proof}
The set $$\{d\in \mathbb{N}\ |\ d\ \text{odd and}\ t_{d}=0\}$$ is syndetic (e.g.~for any $d$ such that $t_d=0$ we have that $t_{4d+3}=0$ and $4d+3$ is odd). 
For each $d$ in the above set, consider $i(d):=2^{n+s} d$; we have $i(d)=2^n \cdot 2^sd$ and $t_{2^sd}=t_d=0$, so $t[i(d),i(d)+2^{n}-1] = t[n]$ and $i(d) = 2^{n+s}d$ takes the desired form (namely, $2^{n+s} \cdot  \text{odd}$).
\end{proof}

We recall the following standard facts about the Thue--Morse sequence. 

\begin{fact} \label{fact}
\begin{enumerate}
    \item For every $n\geq 0$, every factor of $t$ of length $2^{n+2}$ contains an occurrence of the length-$2^n$ prefix $t[n]$. 
    \item Every occurrence of $t[n]$ ($n\geq 1$) in $t$ has starting position $\equiv 0 \mod 2^{n-1}$. Namely, if $t_j\cdots t_{j+2^n-1} = t[n]$ then $j \equiv 0 \mod 2^{n-1}$.
    \item We have $p_t(n) \leq Cn$ for some $C>0$.
    \item $t$ is cube-free.
\end{enumerate}
\end{fact}
For more on the Thue--Morse sequence and related infinite words, see \cite{ASh,ASh2,LothaireBook}.

We define the \emph{interleaved} word
$$ \W = \W_1 \W_2 \W_3 \dots = t_0 w_1 t_1 w_2 t_2 w_3 \cdots \in \{0,1,x_1,\dots,x_b\}^{\mathbb{N}}.$$
For convenience, we index the letters of $\W$ as positive integers.
Namely, $\W_{2i+1} = t_i$ for every $i\geq 0$ and $\W_{2i} = w_i$ for every $i \geq 1$. We note that the idea of interleaving formal languages with  the Thue-Morse sequence has been utilized by Currie and Rampersad in \cite{CR10} to construct infinite binary cube-free words with many squares (that is, an amount exponential in the length).

We let 
\begin{eqnarray*}
\pi\colon \{0,1,x_1,\dots,x_b\}^* & \longrightarrow & \{0,1\}^* \\
\pi'\colon \{0,1,x_1,\dots,x_b\}^* & \longrightarrow & \{x_1,\dots,x_b\}^*
\end{eqnarray*}
be the maps given by erasing $x_1,\dots,x_b$ and $0,1$, respectively; e.g. $$\pi(0x_1 1 x_2 1 x_3 0 x_4)=0110,\ \ \pi'(0x_1 1 x_2 1 x_3 0 x_4)=x_1x_2x_3x_4.$$

\begin{rem} \label{rem:cube-free}
Notice that $\mathcal{W}$ is cube-free, since $t$ is. Indeed, for a subword $\xi$ of $\mathcal{W}$ to have its square also a factor of $\mathcal{W}$, it has to have an even length, so it starts and ends in different alphabets ($\{0,1\}$ vs. $\{x_1,\dots,x_b\}$). Now if $\mathcal{W}[i,j]=\xi^3$ then $\pi(\mathcal{W}[i,j])=\pi(\xi)^3$, contradicting that $t$ is cube-free.
\end{rem}

\begin{prop} \label{prop:interleaved complexity below}
We have $n f(n)\preceq p_\W(n)$.
\end{prop}

\begin{proof}
Fix $n > 1$. We will construct a function
$$ \Psi \colon W(n) \times \{0,1,\dots,2^{n}-1\} \longrightarrow \{\text{Factors of }\W\ \text{of length}\ 2^{n+3}\} $$
and prove that $\Psi$ is at most $2$-to-$1$; consequently, it will follow that $$ p_\W(2^{n+3})\geq \frac{1}{2} \cdot |W(n)| \cdot 2^n \geq 2^{n-1} f(2^{n-1}) $$ (where the last inequality follows from Proposition \ref{prop:define}), and the claim will follow.

Let us now define $\Psi$. For any $u \in W(n)$, by Corollary \ref{cor:every power of 2}, there exists some $m > n$ such that $w_{2^m + 1} \cdots w_{2^m + 2^{n}} = u$ (take the smallest such $m$). This means that $\pi'(\W_{2(2^{m}+1)}\cdots \W_{2(2^m+2^{n})})=u$. For each $0\leq r\leq 2^n - 1$ we let $$\Psi(u,r)=\W_{2(2^m+1-r)}\cdots \W_{2(2^m+1-r)+2^{n+3}-1}.$$

Now fix $u\in W(n)$ and $0\leq r\leq 2^n - 1$ and let $m>n$ be the smallest index for which $w_{2^m + 1} \cdots w_{2^m + 2^{n}} = u$. Given this data, we can compute $\omega:=\Psi(u,r)$ and $\pi'(\omega)=\pi'(\Psi(u,r))$. Suppose that $u'\in W(n)$ and $0\leq r'\leq 2^n - 1$ are such that $\Psi(u,r)=\Psi(u',r')$. Let $m'>n$ be the smallest possible index for which $w_{2^{m'} + 1} \cdots w_{2^{m'} + 2^{n}} = u'$; we aim to show that there is a bounded number of options for such $(u',r')$. 
We have $\Psi(u,r)=\Psi(u',r')$, namely,
\begin{eqnarray*}
\omega & = & \W_{2(2^m+1-r)} \cdots \W_{2(2^m+1-r)+2^{n+3}-1} \\ & = & \W_{2(2^{m'}+1-r')} \cdots \W_{2(2^{m'}+1-r')+2^{n+3}-1}.    
\end{eqnarray*}
Now $\pi(\omega)$ is a subword of $t$ of length $2^{n+2}$. By Fact \ref{fact}, $\pi(\omega)$ contains an occurrence of $t[n]$, say, $\pi(\omega)[j,j+2^{n}-1]=t[n]$. Notice that for every $1\leq k\leq 2^{n+2}$, 
\begin{eqnarray*}
\pi(\omega)[k] & = & \W_{2(2^m+1-r)+2k-1} = t_{2^m-r+k} \\   
\pi(\omega)[k] & = & \W_{2(2^{m'}+1-r')+2k-1} = t_{2^{m'}-r'+k}.
\end{eqnarray*}
It follows that $t[n]$ occurs in $t$ in starting positions
$2^m-r+j$ and $2^{m'}-r'+j$. Again by Fact \ref{fact} it follows that these indices are congruent modulo $2^{n-1}$; since $m,m'>n$, it follows that $r' \equiv r \mod 2^{n-1}$. Hence there are at most $\frac{2^n}{2^{n-1}} = 2$ options for $r'$ (given $r$).
Next, notice that for every $1\leq k\leq 2^{n+2}$, 
\begin{eqnarray*}
\pi'(\omega)[k] & = & \W_{2(2^m+1-r)+2k-2} = w_{2^m-r+k} \\   
\pi'(\omega)[k] & = & \W_{2(2^{m'}+1-r')+2k-2} = w_{2^{m'}-r'+k}.
\end{eqnarray*}
In particular, $\pi'(\omega)[r'+1,r'+2^n]=u'$. Since we have at most $2$ options for $r'$, and $\pi'(\omega)=\pi'(\Psi(u,r))$ can be computed given $(u,r)$, there are now at most $2$ options for $u'$ (one for each option for $r'$; notice that if $r'=r$ then $u'=u$). It follows that $\Psi$ is at most $2$-to-$1$, and the claim follows.
\end{proof}

\begin{prop} \label{prop:interleaved complexity above}
We have $p_\W(n)\preceq n f(n)$.
\end{prop}

\begin{proof}
For every factor $\omega$ of $\W$ of length $2^{n+3}$, we have that $\pi(\omega)$ (a factor of $t$ of length $2^{n+2}$) contains an occurrence of $t[n]$, and all occurrences of $t[n]$ within it have starting positions \emph{in $t$} that are congruent to $0$ modulo $2^{n-1}$ (by Fact \ref{fact}), corresponding to starting positions in $\W$ that are congruent to $1$ modulo $2^n$ (recall that $t_i = \W_{2i+1}$). We thus have at most $\frac{2^{n+3}}{2^n}=8$ candidates for the starting position of $\omega$ in $\W$ modulo $2^{n+3}$. By Lemma \ref{lem:starting position modulo}, each of these possible starting positions gives us at most $4f(2^{n+3})$ options for $\pi'(\omega)$; and we have at most $C\cdot 2^n$ options for $\pi(\omega)$ (by Fact \ref{fact}). Altogether, we have at most $$ C\cdot 2^n \cdot 8 \cdot 4f(2^{n+3}) $$ possible factors of $\W$ of length $2^{n+3}$. It follows that $p_\W(n)\preceq n f(n)$.
\end{proof}

\begin{prop} \label{prop:interleaved minimal}
The word $\W$ is uniformly recurrent.
\end{prop}
\begin{proof}
It suffices to prove that for every $k$, the set of indices in which the prefix $\W[1,2^{k+1}]=t_0 w_1 t_1 w_2 \cdots t_{2^k-1} w_{2^k}$ appears in $\W$ is syndetic. By Corollary \ref{cor:every power of 2}, there exists some $l>k$ such that for every odd $d\geq 1$ we have $w_{d 2^l+1}\dots w_{d 2^l+2^k}=w_1 \dots w_{2^k}$.

Recall that, by Lemma \ref{lem:syndeticTM} for $n=k,\ s=l-k$, the set 
$$\mathcal{S} := \{i\in \mathbb{Z}_{\geq 0}\ |\ t[i,i+2^k-1]=t_0 t_1\dots t_{2^k-1}\ \text{and}\ i=2^l\cdot \text{odd} \}$$
is syndetic. Notice that for each $i\in \mathcal{S}$ we have that $i+1$ is of the form $d 2^l + 1$ for odd $d$, so $w_{d 2^l + 1}\dots w_{d 2^l + 2^k}=w_1\dots w_{2^k}$.
Now 
\begin{align*}
\W_{2i+1}\W_{2i+2}\dots \W_{2i+2^{k+1}} & =t_i w_{i+1} t_{i+1} w_{i+2} \dots t_{i+2^k-1} w_{i+2^k} \\ & =t_0 w_1 t_1 w_2 \dots t_{2^k-1} w_{2^k}
\end{align*}
so the set $$\{j\in \mathbb{N}\ |\ \W[j,j+2^{k+1}-1]=\W[1,2^{k+1}] \}$$ contains $2\mathcal{S}+1$ and is thus syndetic itself.
\end{proof}

\subsection{Unique ergodicity}

We now aim to prove that $\W$ generates a uniquely ergodic subshift. Toward this end, we must prove the following. Let $u=u_1\dots u_k,v=v_1\dots v_l$ be any words over the alphabet $\{0,1,x_1,\dots,x_b\}$. Denote  $\Phi_u(v)=|\{1\leq i\leq l-k+1|v_i\dots v_{i+k-1}=u\}|$, the number of occurrences of $u$ in $v$, and $\varphi_u(v)=\frac{\Phi_u(v)}{l} \in [0,1]$. We must prove that for every $u$ there exists some constant $\chi_u$ measuring the `frequency' of $u$ within $\W$, namely, for any sequence of factors $v_1,v_2,\dots$ of $\W$ with $|v_i|\xrightarrow{i\rightarrow \infty} \infty$, we have $\varphi_u(v_i)\xrightarrow{i\rightarrow \infty} \chi_u$. This property is equivalent to having the subshift generated by $\W$ uniquely ergodic.

Let $m\in \mathbb{N},\ \alpha=\alpha_1\dots \alpha_{2^m}\in \Sigma^{2^m},\ \xi=d_0\dots d_{2^m-1}\in \{0,1\}^{2^m}$ be given. We denote $\alpha \star \xi=d_0\alpha_1\dots d_{2^m-1}\alpha_{2^m}$. For $v=v_1\dots v_{2^k}\in \Sigma^{2^k},\ k>m$, we denote
$$N_v(\alpha,\xi)=|\{0\leq j\leq 2^{k}-2^m\ |\ t_j v_{j+1} t_{j+1}\dots t_{j+2^m-1} v_{j+2^m}=\alpha \star \xi\}|.$$
Here, as in the previous subsections, $t=t_0t_1\dots$ is the Thue-Morse sequence. Notice that $$N_v(\alpha,\xi)=\Phi_{\alpha \star \xi}(v \star t[k])\ \ \text{and}\ \ N_v(\alpha,\bar{\xi})=\Phi_{\alpha\star \bar{\xi}}(v\star t[k])=\Phi_{\alpha\star \xi}(v\star \overline{t[k]}).$$

Let
\begin{eqnarray*}
m_k(\alpha,\xi) & = & \min\{N_v(\alpha,\xi)\ |\ v\in W(k)\} \\
\tilde{m}_k(\alpha,\xi) & = & \min\{m_k(\alpha,\xi),m_k(\alpha,\bar{\xi})\} \\
M_k(\alpha,\xi) & = & \max\{N_v(\alpha,\xi)\ |\ v\in W(k)\} \\
\tilde{M}_k(\alpha,\xi) & = & \max\{M_k(\alpha,\xi),M_k(\alpha,\bar{\xi})\} \\
A_k(\alpha,\xi) & = & \frac{\tilde{m}_k(\alpha,\xi)}{2^k} \\ 
B_k(\alpha,\xi) & = & \frac{\tilde{M}_k(\alpha,\xi)+2^m-1}{2^k}.
\end{eqnarray*}

Notice that $ A_k(\alpha,\xi) \leq B_k(\alpha,\xi) $ for all $k>m$.

\begin{lem}
We have $A_k(\alpha,\xi)\leq A_{k+1}(\alpha,\xi)$ and $B_k(\alpha,\xi)\geq B_{k+1}(\alpha,\xi)$.
\end{lem}

\begin{proof}
Suppose $v\in W(k+1)$ is given. Since $W(k+1)=W(k)C(k)\subseteq W(k)W(k)$, we can decompose $v=v' v''$ for $v',v''\in W(k)$. Now $v\star t[k+1]=(v'\star t[k])(v'' \star \overline{t[k]})$ so $$N_v(\alpha,\xi)\geq N_{v'}(\alpha,\xi)+N_{v''}(\alpha,\bar{\xi})\geq 2\tilde{m}_k(\alpha,\xi).$$ Likewise, $$N_v(\alpha,\bar{\xi})\geq N_{v'}(\alpha,\bar{\xi})+N_{v''}(\alpha,\xi)\geq 2\tilde{m}_k(\alpha,\xi).$$ Hence 
\begin{eqnarray*} A_{k+1}(\alpha,\xi) & = & 2^{-(k+1)}\cdot \min\{N_v(\alpha,\xi),N_v(\alpha,\bar{\xi})\ |\ v\in W(k+1)\} \\ & \geq &  2^{-(k+1)} \cdot  2\tilde{m}_k(\alpha,\xi) = A_k(\alpha,\xi).\end{eqnarray*}

Next, given any $v\in W(k+1)$ such that $\tilde{M}_{k+1}(\alpha,\xi)=N_v(\alpha,\xi)$ or $\tilde{M}_{k+1}(\alpha,\xi)=N_v(\alpha,\bar{\xi})$, decompose $v=v' v''$ as before. Notice that $$N_v(\alpha,\xi)\leq N_{v'}(\alpha,\xi)+N_{v''}(\alpha,\bar{\xi})+2^m-1$$ (any occurrence of $\alpha\star \xi$ in $v \star t[k+1]$ is either within $v'\star t[k],v''\star \overline{t[k]}$, or has its starting position in the $2^m-1$ `stitch' letters in the end of $v'$). Hence $$N_v(\alpha,\xi)+2^m-1\leq 2\left(\tilde{M}_k (\alpha,\xi)+2^m-1\right);$$ the same inequality holds for $N_v(\alpha,\bar{\xi})$, so $$\tilde{M}_{k+1}(\alpha,\xi)+2^m-1 \leq 2\left(\tilde{M}_k (\alpha,\xi)+2^m-1\right),$$ and thus $B_{k+1}(\alpha,\xi)\leq B_k(\alpha,\xi)$.
\end{proof}

It follows that $\{A_k(\alpha,\xi)\}_{k=m+1}^{\infty},\{B_k(\alpha,\xi)\}_{k=m+1}^{\infty}$ are both convergent.

\begin{prop} \label{prop:A=B}
We have $\lim_{k\rightarrow \infty} A_k(\alpha,\xi)=\lim_{k\rightarrow \infty} B_k(\alpha,\xi)$.
\end{prop}

\begin{proof}
Recall that by the construction of $\{W(k),C(k)\}_{k=1}^{\infty}$, we have that $C(k)$ is a singleton infinitely often. Fix a sequence $m<j_1<j_2<\dots$ such that $j_{i+1}-j_i\geq 2$ and $C(j_i)=\{v_i\}$ for some $v_i\in W(j_i)$. 
Thus every $v\in W(j_i+2)$ decomposes as $v=v'v_iv''v_i$ for some $v',v''\in W(j_i)$. 
Recall that the Thue-Morse sequence enjoys the property (in fact, definition) $t[n+1]=t[n]\overline{t[n]}$. 
Hence for every $v\in W(j_i+2)$, we have that 
\begin{align} \label{eq:decompose}
v \star t[j_i+2] & =  v' v_i v'' v_i \star t[j_i]\overline{t[j_i]}\overline{t[j_i]}t[j_i] \nonumber \\ & =  \left(v' \star t[j_i]\right) \left(v_i \star \overline{t[j_i]}\right) \left(v''\star \overline{t[j_i]}\right) \left(v_i\star t [j_i]\right)
\end{align} (for some $v',v''\in W(j_i)$) contains disjoint occurrences of $v_i\star t[j_i]$, in which we count $N_{v_i}(\alpha,\xi)$ occurrences of $\alpha\star \xi$, and of $v_i\star \overline{t[j_i]}$, in which we count $\Phi_{\alpha\star\xi}(v_i \star \overline{t[j_i]})=N_{v_i}(\alpha,\bar{\xi})$ occurrences of $\alpha\star \xi$. In addition, we have two other (all are pairwise disjoint) subwords, $v'\star t[j_i]$ and $v''\star \overline{t[j_i]}$, each of which contributes at least $\tilde{m}_{j_i}(\alpha,\xi)$. Therefore
\begin{equation} \label{eq:mji} m_{j_i+2}(\alpha,\xi)\geq 2\tilde{m}_{j_i}(\alpha,\xi)+N_{v_i}(\alpha,\xi)+N_{v_i}(\alpha,\bar{\xi}).
\end{equation}
Moreover, again by (\ref{eq:decompose}), we see that any occurrence of $\alpha\star \xi$ in $v\star t[j_i+2]$ is either inscribed in one of $v' \star t[j_i],v_i \star \overline{t[j_i]},v''\star \overline{t[j_i]},v_i\star t [j_i]$, or its starting position is in one of the three `stitch' $2^m-1$ positions in the end of each one of the first three of those. We have at most $\tilde{M}_{j_i}(\alpha,\xi),N_{v_i}(\alpha,\bar{\xi}),\tilde{M}_{j_i}(\alpha,\xi),N_{v_i}(\alpha,\xi)$ occurrences of $\alpha\star \xi$ in $v' \star t[j_i],v_i \star \overline{t[j_i]},v''\star \overline{t[j_i]},v_i\star t [j_i]$, respectively, and therefore altogether
\begin{equation} \label{eq:Mji}
M_{j_i+2}(\alpha,\xi)\leq 2\tilde{M}_{j_i}(\alpha,\xi)+N_{v_i}(\alpha,\xi)+N_{v_i}(\alpha,\bar{\xi})+3(2^m-1).
\end{equation}
By analogous arguments, we also get
\begin{equation} \label{eq:mbji} m_{j_i+2}(\alpha,\bar{\xi})\geq 2\tilde{m}_{j_i}(\alpha,\xi)+N_{v_i}(\alpha,\xi)+N_{v_i}(\alpha,\bar{\xi})
\end{equation}
and
\begin{equation} \label{eq:Mbji}
M_{j_i+2}(\alpha,\bar{\xi})\leq 2\tilde{M}_{j_i}(\alpha,\xi)+N_{v_i}(\alpha,\xi)+N_{v_i}(\alpha,\bar{\xi})+3(2^m-1).
\end{equation}
Hence, by (\ref{eq:mji}) and (\ref{eq:mbji}),
\begin{equation} \label{eq:mtji}
\tilde{m}_{j_i+2}(\alpha,\xi)\geq 2\tilde{m}_{j_i}(\alpha,\xi)+N_{v_i}(\alpha,\xi)+N_{v_i}(\alpha,\bar{\xi})
\end{equation}
and by (\ref{eq:Mji}) and (\ref{eq:Mbji}),
\begin{equation} \label{eq:Mtji}
\tilde{M}_{j_i+2}(\alpha,\bar{\xi})\leq 2\tilde{M}_{j_i}(\alpha,\xi)+N_{v_i}(\alpha,\xi)+N_{v_i}(\alpha,\bar{\xi})+3(2^m-1).
\end{equation}
Define $$\lambda_i = \frac{N_{v_i}(\alpha,\xi)+N_{v_i}(\alpha,\bar{\xi})}{2^{j_i+2}}.$$
By (\ref{eq:mtji}),
\begin{eqnarray*}
A_{j_i+2}(\alpha,\xi) & = & \frac{\tilde{m}_{j_i+2}(\alpha,\xi)}{2^{j_i+2}} \geq  \frac{2\tilde{m}_{j_i}(\alpha,\xi)+N_{v_i}(\alpha,\xi)+N_{v_i}(\alpha,\bar{\xi})}{2^{j_i+2}} \\ & = & 
\frac{A_{j_i}(\alpha,\xi)}{2} + \lambda_i
\end{eqnarray*}
and by (\ref{eq:Mtji}),
\begin{eqnarray*}
B_{j_i+2}(\alpha,\xi) & = & \frac{\tilde{M}_{j_i+2}(\alpha,\xi)+2^m-1}{2^{j_i+2}} \\ & \leq & \frac{2\tilde{M}_{j_i}(\alpha,\xi)+N_{v_i}(\alpha,\xi)+N_{v_i}(\alpha,\bar{\xi})+4(2^m-1)}{2^{j_i+2}} \\ & = & \frac{B_{j_i}(\alpha,\xi)}{2} + \lambda_i + \frac{2^m-1}{2^{j_i + 1}} \\ & < & \frac{B_{j_i}(\alpha,\xi)}{2} + \lambda_i +2^{m-j_i}  \leq \frac{B_{j_i}(\alpha,\xi)}{2} + \lambda_i + 2^{-i}.
\end{eqnarray*}
Hence
\begin{eqnarray*}
A_{j_{i+1}}(\alpha,\xi) &  \geq & A_{j_i+2}(\alpha,\xi) \geq \frac{A_{j_i}(\alpha,\xi)}{2} + \lambda_i \\
B_{j_{i+1}}(\alpha,\xi) & \leq & B_{j_i+2}(\alpha,\xi) \leq \frac{B_{j_i}(\alpha,\xi)}{2}+\lambda_i+2^{-i}.
\end{eqnarray*}
Denote $\delta_i = B_{j_i} - A_{j_i}$, and observe that $0\leq \delta_i\leq 1$. 
Thus
$$\delta_{i+1} = B_{j_{i+1}} - A_{j_{i+1}} \leq \frac{B_{j_i}(\alpha,\xi)-A_{j_i}(\alpha,\xi)}{2} + 2^{-i} = \frac{\delta_i}{2} + 2^{-i},$$
so $\delta_i \xrightarrow{i\rightarrow \infty} 0$ and $\lim_{k\rightarrow \infty} A_k(\alpha,\xi) = \lim_{k\rightarrow \infty} B_k(\alpha,\xi)$, as claimed.
\end{proof}

\begin{prop} \label{prop:forstar}
Let $\alpha\in \Sigma^{2^m},\xi\in \{0,1\}^{2^m}$ be given, and suppose that $\alpha \star \xi$ factors $\W$.  There exists a constant $\rho=\rho_{\alpha,\xi}\in (0,1]$ such that
$$ \frac{1}{|U|} \Phi_{\alpha\star \xi}(U) \xrightarrow{|U|\rightarrow \infty} \rho, $$
regardless of which $U$ are chosen, as long as they are factors of $\W$ and their lengths tend to infinity.
\end{prop}
\begin{proof}
Given $\alpha\in \Sigma^{2^m},\xi\in \{0,1\}^{2^m}$, denote $$\chi_{\alpha \star \xi} := \lim_{k\rightarrow \infty} A_k(\alpha,\xi) = \lim_{k\rightarrow \infty} B_k(\alpha,\xi)$$
(utilizing Proposition \ref{prop:A=B}.)
This means that, for every $\varepsilon>0$ (assume $\varepsilon<\chi_{\alpha\star \xi}$), there exists $k_0$ such that for all $k\geq k_0$ and for all $v\in W(k)$, $$\Big|\Phi_{\alpha\star \xi}(v\star t[k])-2^k \chi_{\alpha\star \xi}\Big|,\Big|\Phi_{\alpha\star \xi}(v\star \overline{t[k]})-2^k \chi_{\alpha\star \xi}\Big|<2^k \varepsilon.$$
Further assume that $k_0>m+\log_2 \frac{1}{\varepsilon}$ and $k_0>2$. Suppose $L>\frac{2^{k_0+10}}{\varepsilon}$ and let $U$ be an arbitrary factor of $$ \W = \W_1 \W_2 \W_3 \dots = t_0 w_1 t_1 w_2 t_2 w_3 \cdots$$ of length $L$.
Write 
\begin{eqnarray*} U =  U_0 & & \W[D\cdot 2^{k_0+1} + 1,(D+1)\cdot 2^{k_0+1}] \\ & & \W[(D+1)\cdot 2^{k_0+1} + 1,(D+2)\cdot 2^{k_0+1}] \\ & & \vdots \\ & & \W[(D+r)\cdot 2^{k_0+1} + 1,(D+r+1)\cdot 2^{k_0+1}] \ \ \ U_1
\end{eqnarray*}
for some $U_0,U_1$ of lengths at most $2^{k_0+1}$ and for some $D,r\geq 0$. 
Thus 
\begin{eqnarray*}
(r+1)\cdot 2^{k_0+1} \leq L & = & |U| \\ & = & |U_0|+|U_1|+(r+1)\cdot 2^{k_0+1} \\  & \leq & (r+3)\cdot 2^{k_0+1}.
\end{eqnarray*}
Notice that each one of the above blocks between $U_0,U_1$ takes the form
\begin{eqnarray*}
& & \W[(D+l)\cdot 2^{k_0+1} + 1,(D+l+1)\cdot 2^{k_0+1}] = \\ & & t_{(D+l)\cdot 2^{k_0}} w_{(D+l)\cdot 2^{k_0}+1} t_{(D+l)\cdot 2^{k_0}+1}\dots w_{(D+l+1)\cdot 2^{k_0}} = \\ & & \underbrace{w_{(D+l)\cdot 2^{k_0}+1}\dots w_{(D+l+1)\cdot 2^{k_0}}}_{\in W(k_0)} \star \underbrace{t_{(D+l)\cdot 2^{k_0}}\dots t_{(D+l+1)\cdot 2^{k_0}-1}}_{\in \{t[k_0],\overline{t[k_0]}\}}
\end{eqnarray*}
for some $0\leq l\leq r$. (Notice that, indeed, $t[K\cdot 2^{k_0},(K+1)\cdot 2^{k_0}-1]\in \{ t[k_0],\overline{t[k_0]} \}$ by the inductive definition of the Thue-Morse sequence.) Hence $$\Big| \Phi_{\alpha\star \xi}(\W[(D+l)\cdot 2^{k_0+1} + 1,(D+l+1)\cdot 2^{k_0+1}])-2^{k_0}\chi_{\alpha\star \xi} \Big| < 2^{k_0} \varepsilon.$$ 
Therefore
$$\Phi_{\alpha\star \xi}(U) \geq (r+1) \cdot \left(2^{k_0}\chi_{\alpha\star \xi} - 2^{k_0} \varepsilon\right) $$
and
\begin{eqnarray*}
\frac{1}{L}\Phi_{\alpha \star \xi}(U) & \geq & \frac{(r+1) \cdot \left(2^{k_0}\chi_{\alpha\star \xi} - 2^{k_0} \varepsilon \right)}{(r+3)\cdot 2^{k_0+1}} \\ & \geq & (1-\varepsilon)\left(\frac{\chi_{\alpha\star \xi}}{2} - \frac{\varepsilon}{2}\right) \geq (1-\varepsilon)\cdot \frac{\chi_{\alpha \star \xi}}{2} - \varepsilon.
\end{eqnarray*}
Now
\begin{eqnarray*}
\Phi_{\alpha\star \xi}(U) & \leq & |U_0|+|U_1|+(r+1)\cdot \left(2^{k_0} \chi_{\alpha\star \xi} + 2^{k_0} \varepsilon\right) + \underbrace{(r+1)\cdot 2^m}_{\substack{\text{Starting positions} \\ \text{in the end of} \\ \text{each block}}} \\
& \leq & 2^{k_0+2} + (r+1)\cdot \left(2^{k_0} \chi_{\alpha\star \xi} + 2^{k_0} \varepsilon\right) + (r+1) \cdot 2^m \\ 
& \leq & (r+1)\cdot \left(2^{k_0} \chi_{\alpha\star \xi} + 2^{k_0}  \varepsilon\right) + 2^{k_0+2} + (r+1)\cdot \varepsilon \cdot 2^{k_0} \\ & \leq & (r+1)\cdot \left(2^{k_0} \chi_{\alpha\star \xi} + 2^{k_0} \varepsilon\right) + \frac{\varepsilon L}{256} + \frac{\varepsilon L}{2}. 
\end{eqnarray*}
Therefore
$$ \frac{1}{L}\Phi_{\alpha\star \xi}(U) \leq \frac{(r+1)\cdot \left(2^{k_0} \chi_{\alpha\star \xi} + 2^{k_0} \varepsilon\right)}{(r+1)\cdot 2^{k_0+1}} + \frac{3\varepsilon}{4} = \frac{\chi_{\alpha\star \xi}}{2} + \frac{5 \varepsilon}{4}. $$
For conclusion, we proved that 
$$ \frac{1}{|U|} \Phi_{\alpha\star \xi}(U) \xrightarrow {|U|\rightarrow \infty} \frac{\chi_{\alpha\star \xi}}{2}, $$
regardless of the choice of $U$ (as long as their lengths tend to infinity). We know that $\frac{1}{|U|}\Phi_{\alpha\star \xi}(U)$ is bounded away from zero, since $\alpha\star \xi$ factors $\W$, and since $\W$ is uniformly recurrent (Proposition \ref{prop:interleaved minimal}).
\end{proof}

\begin{prop} \label{prop:frequencies}
Let $u$ be any finite factor of $\W$. There exists a constant $\rho=\rho_u\in (0,1]$ such that
$$ \frac{1}{|U|} \Phi_{u}(U) \xrightarrow{|U|\rightarrow \infty} \rho, $$
regardless of which $U$ are chosen, as long as they are factors of $\W$ and their lengths tend to infinity.
\end{prop}

\begin{proof}
By Proposition \ref{prop:forstar}, we know that the claim holds for $u$ taking the form $\alpha\star \xi$ for $\alpha,\xi$ of length a power of $2$.
Let $u$ be any finite factor of $\W$; it factors some word of the form $\alpha \star \xi$ with $\alpha\in \Sigma^{2^m},\xi\in \{0,1\}^{2^m}$. Write $\alpha \star \xi = vuv'$ and denote $|v|=p,|v'|=q$, so $\alpha\star \xi[p+1,p+|u|]=u$. 
Let $$S = \{\alpha \star \xi\ \text{factor of}\ \W\ |\ \alpha \in \Sigma^{2^m},\xi \in \{0,1\}^{2^m},\ \alpha \star \xi[p+1,p+|u|]=u\}$$
and observe that $S\neq \emptyset$. Then $\Phi_u(U) \geq \sum_{\alpha\star \xi\in S} \Phi_{\alpha \star \xi}(U)$, as any occurrence of $\alpha \star \xi\in S$ within $U$ contains a (unique) occurrence of $u$ starting from its $p+1$-th letter, and no distinct words from $S$ can have the same occurrence of $u$ starting in their $p+1$-th position (or they must coincide). 

Conversely, any occurrence of $u$ in $U$ starting in the $i$-th letter of $U$, namely, $U[i,i+|u|-1]=u$, can be (uniquely) prolonged by $p$ to the left and by $q$ to the right to a word $\alpha\star \xi\in S$, if $i\geq p+1$ and $i+|u|-1\leq |U|-q$. The remaining cases are $1\leq i\leq p$ and $|U|-|u|-q+1 < i \leq |U|-|u|+1$ (notice that the starting position of $u$ in $U$ cannot exceed $|U|-|u|+1$). Hence
$$\Phi_u(U)\leq \sum_{\alpha \star \xi\in S} \Phi_{\alpha \star \xi}(U) + p + q,$$ so
$$\frac{1}{|U|} \sum_{\alpha \star \xi \in S} \Phi_{\alpha \star \xi}(U) \leq \frac{1}{|U|} \Phi_u(U) \leq \frac{1}{|U|} \sum_{\alpha \star \xi\in S} \Phi_{\alpha \star \xi}(U) + \frac{p + q}{|U|}$$
and, subsequently,
$$\frac{1}{|U|}\Phi_u(U) \xrightarrow{|U|\rightarrow \infty} \sum_{\alpha \star \xi \in S} \rho_{\alpha \star \xi}$$
and the claim follows (the limit is positive because $S$ contains at least one factor of $\W$ extending $u$).
\end{proof}

\begin{proof}[{Proof of Theorem \ref{thm:ergodic}}]
Given a non-decreasing, submultiplicative function $f\colon \mathbb{N}\rightarrow \mathbb{N}$, consider the above construction $\W$ (notice that we may assume that $f$ is unbounded and subexponential; indeed, if $f$ is bounded then $nf(n)\sim n$, and if $f(n)$ is exponential, so is $nf(n)$, and subshifts as in the theorem's assertion of linear and of exponential complexity are well known). By Propositions \ref{prop:interleaved complexity below} and \ref{prop:interleaved complexity above} we have that $p_\W(n) \sim n f(n)$. By Proposition \ref{prop:interleaved minimal}, the word $\W$ is uniformly recurrent, so it generates a minimal subshift. By Proposition \ref{prop:frequencies}, factors have well-defined frequencies within $\W$, and thus $\W$ generates a uniquely ergodic subshift.  By Remark \ref{rem:cube-free}, it is also cube-free.
\end{proof}

\subsection{Applications to the growth of algebraic structures}

Let $S=\langle X \rangle$ be a finitely generated semigroup with a finite generating set $X$. The growth function of $S$ with respect to $X$ is given by $$\gamma_{S,X}(n)=|X^{\leq n}|=|\{\text{Products of at most}\ n\ \text{generators}\}|.$$
When $S$ is a group, one takes $X$ to be a symmetric generating set to obtain the standard notion of the growth function of a group. The first groups with super-polynomial, subexponential growth were constructed by Grigorchuk \cite{Gri1}.

Let $K$ be a field and let $A$ be a finitely generated $K$-algebra with a finite generating set $X$. The growth of $A$ with respect to $X$ is given by $$\gamma_{A,X}(n)=\dim_K X^{\leq n} = \dim_K \Span\{\text{Products of at most}\ n\ \text{generators}\}.$$

In all of the above cases, the growth function is independent of the choice of a generating set when considered up to bi-Lipschitz equivalence. For more on the growth of algebras, see \cite{KL}.

The classes of growth functions of semigroups and associative algebras coincide with the class of growth functions of hereditary languages; these classes have been fully characterized by Bell and the fourth-named author in \cite{BZ}. These growth functions are always non-decreasing and submultiplicative (assuming a neutral element). Bergman's gap theorem \cite{Ber} asserts that the growth function of a semigroup (equivalently, an associative algebra) is either at least quadratic, or linear, or bounded. The following can thus be viewed as a constructive converse to Bergman's gap theorem:

\begin{cor} \label{cor:simple2}
    Let $f\colon \mathbb{N}\rightarrow \mathbb{N}$ be a non-decreasing, submultiplicative function. Then there exists a finitely generated simple associative algebra $A$ such that $\gamma_A(n)\sim n^2 f(n)$.
\end{cor}

\begin{proof}
By Theorem \ref{thm:ergodic}, there exists an infinite, uniformly recurrent word $w$ over a finite alphabet such that $p_w(n)\sim nf(n)$. We can replace $w$ by a bi-infinite word. 
Nekrashevych \cite{Nek} associated with every aperiodic uniformly recurrent bi-infinite word $w$, a finitely generated simple algebra $A_w$ (namely, the convolution algebra of the groupoid of the action of $\mathbb{Z}\curvearrowright \overline{\text{shift orbit of}\ w}$) such that $\gamma_{A_w}(n) \sim np_{w}(n)$. Thus $\gamma_{A_w}(n) \sim np_w(n)\sim n^2 f(n)$.
\end{proof}

Notice that this strengthens \cite[Corollary~5.5]{GMZ}, providing a precise realization of $n^2 f(n)$ as the growth function of a simple algebra rather than an approximation.
We conclude with the following application to the growth of semigroups satisfying semigroup laws. 
Shneerson studied the possible growth rates of semigroups satisfying semigroup laws, see e.g. \cite{Shn04}, the comprehensive paper \cite{Shn08}, and references therein; in \cite{BR08}, an interesting example of a semigroup of polynomial growth of irrational degree, satisfying the law $X^6=X^4$, is given.

\begin{cor}
    Let $f\colon \mathbb{N}\rightarrow \mathbb{N}$ be a non-decreasing, submultiplicative function. Then there exists a finitely generated just-infinite semigroup with zero $S$ satisfying the law $X^3=0$, of growth $\gamma_S(n)\sim n^2 f(n)$.
\end{cor}
\begin{proof}
By Theorem \ref{thm:ergodic}, there exists an infinite, uniformly recurrent word $w$ over a finite alphabet $\Sigma$ such that $p_w(n)\sim nf(n)$.
Consider the monomial semigroup $S=\langle \Sigma\ |\ u=0\ \text{if}\ u\ \text{is not a subword of}\ w\rangle$. Now 
\begin{eqnarray*}
\gamma_{S,X}(n) & = & |\{\text{Subwords of}\ w\ \text{of length}\leq n\}| \\ & = & \sum_{k=0}^n p_w(k) \sim np_w(n)\sim n^2 f(n).
\end{eqnarray*}
That $S$ satisfies the semigroup law $X^3=0$ is evident from the fact that $w$ is cube-free. That $S$ is just-infinite follows from the uniform recurrence of $w$ by means of standard reduction processes (see \cite{BBL} for more on monomial semigroups and algebras), which we include here for completeness. Fix an arbitrary order on $\Sigma$ and extend it to the degree-lexicographic order (denoted $\prec$) on the set of all words in $\Sigma$. Let $\varphi \colon S\to T$ be a proper homomorphic image and let $\varphi(u)=\varphi(v)$ be a relation in $T$ that does not hold in $S$, namely, $u\neq v$. We may assume that $u\neq 0$, and that $v\prec u$. Hence $u$ is a subword of $w$, and, by uniform recurrence, there exists some $N$ such that every subword of $w$ of length $N$ contains an occurrence of $u$. We claim that $\Image(\varphi)=\varphi\left(\{\text{Words of length} < N\}\right) \cup \{0\}$. Let $t$ be a subword of $w$ of length $\geq N$ such that $\varphi(t)\neq 0$ and $\varphi(t)\neq \varphi(t')$ for any $|t'| < N$. Consider such $t$ minimal with respect to $\prec$ (this is a well-order). Now $t$ must contain an occurrence of $u$, say, $t=puq$, which means that $\varphi(t)=\varphi(pvq)$. Since $pvq \prec puq=t$, it follows that $pvq$ is not a subword of $w$, or else, this contradicts the minimality of $t$; therefore $\varphi(pvq)=0$, so $\varphi(t)=0$. It follows that $T=\Image(\varphi)$ is finite.
\end{proof}

For context, there is no known finitely generated group of intermediate growth (that is, subexponential and super-polynomial) that satisfies a non-trivial group law.

\end{document}